%% file: main.tex
\documentclass[11pt,a4paper,twoside]{amsart}
\usepackage[top=1in, bottom=1.25in, marginparwidth=1in, inner=1in, outer=1in]{geometry}

\usepackage{lmodern}
\usepackage[T2A,T1]{fontenc} 
\usepackage[utf8]{inputenc}
\usepackage{inputenc}
\usepackage[UKenglish]{babel}
\usepackage{yfonts}
\usepackage{bm}
\usepackage{booktabs}
\usepackage{afterpage}
\usepackage{pinlabel}
\usepackage{ifthen}
\usepackage{tabularx}
\usepackage{nicefrac}

 \usepackage{picins}

\makeatletter
\newcommand*{\math@version@bold}{bold}
\DeclareMathOperator\DD{
	\textrm{%
		\usefont{T2A}{cmr}{\ifx\math@version\math@version@bold bx\else m\fi}{n}%
		\CYRD
	}%
} 
\makeatother

\usepackage[labelfont=bf]{caption}
\DeclareCaptionFormat{myformat}{\fontsize{9}{10}\selectfont#1#2#3}
\usepackage{soul}
\usepackage{calc}
\usepackage{float}
\usepackage{graphicx}
\usepackage[table]{xcolor}
\colorlet{lightred}{red!20!white}
\colorlet{lightblue}{blue!30!white}
\colorlet{darkgreen}{green!70!black}
\colorlet{lightgreen}{green!50!white}
\colorlet{gold}{yellow!90!black!70!red}

\usepackage{multicol,tabu,longtable}
\usepackage{rotating}
\usepackage{tikz}
\usepackage{tikz-cd}
\usetikzlibrary{shapes,decorations}
\tikzcdset{arrow style=tikz, diagrams={>=latex}}
\tikzset{>=latex}

\usepackage{caption}
\usepackage[labelfont=up,margin=5pt]{subcaption}
\usepackage{wrapfig}
\usepackage[section]{placeins}
\usepackage{setspace} 
\usepackage{hyperref}
\hypersetup{%
	pdfencoding=auto, %
	bookmarksdepth=3, %
	colorlinks=true, %
	citecolor=black, %
	urlcolor=black, %
	linkcolor=black %
}
\usepackage{bookmark}

\usepackage{amsmath}
\usepackage{amsfonts}
\usepackage{amssymb}
\usepackage{amsthm}
\usepackage{mathtools}
\usepackage{mathrsfs}
\mathtoolsset{mathic=true}

\newtheorem{theorem}{Theorem}[section]
\newtheorem*{theorem*}{Theorem}
\newtheorem{lemma}[theorem]{Lemma}
\newtheorem{conjecture}[theorem]{Conjecture}
\newtheorem{question}[theorem]{Question}

\newtheorem{questions}[theorem]{Questions}

\theoremstyle{definition}
\newtheorem{definition}[theorem]{Definition}
\newtheorem{example}[theorem]{Example}
\newtheorem{counterexample}[theorem]{Counterexample}

\newtheorem{remark}[theorem]{Remark}

\usepackage[makeroom]{cancel}
\usepackage{enumitem}
\newcounter{dummy}
\makeatletter
\newcommand\myitem[1][]{%
	\item[\textnormal{(}#1\textnormal{)}]\refstepcounter{dummy}%
	\def\@currentlabel{\textnormal{(}#1\textnormal{)}}%
}
\makeatother

\usepackage[normalem]{ulem} 

\input{notation_commands.tex}

\input{PSTricks/input}

\newcommand{\mynode}[3]{%
	\ifthenelse{\equal{#3}{1}}{%
		\path [fill=lightgray,draw=none] (#1,#2) rectangle ++(-1,-1);
	}{%
		\path [fill=lightgray,draw=none] (#1,#2) rectangle ++(-1,-1) node [midway]{\tiny\( #3 \)};%
	}
}
\newcommand{\mygrid}[4]{%
	\draw[gray,step=1] (#1-1,#2-1) grid (#3,#4);
	\draw[gray] (#1-1,#2-1) -- (#1-1.25,#2-1);
	\draw[gray] (#1-1,#2-1) -- (#1-1,#2-1.25);
	\draw (#1-1,#2-1) node[left] {\tiny\( #2 \)};
	\path [draw=none] (#1-0.5,#2-1.25) rectangle ++(-1,-1) node [midway]{\tiny\( #1 \)};
}
\newcommand{\mylabel}[5]{%
	\draw (#1-6,#2-#2/2+#4/2-0.5) node[right] {\( #5 \)};
}

\newcommand{\myrow}[5]{%
\begin{tabularx}{73pt}{%
		>{\centering\arraybackslash}m{14pt}
		>{\centering\arraybackslash}m{19pt}%
		>{\centering\arraybackslash}m{10pt}}%
	#1 \(#2\)\!\!\! & \(#3\) & \(#4\) #5
\end{tabularx} 
}

\begin{document}
\title{Khovanov homology and strong inversions}

\author{Artem Kotelskiy}
\address{Department of Mathematics \\ Indiana University}
\email{artofkot@iu.edu}

\author{Liam Watson}
\address{Department of Mathematics \\ University of British Columbia}
\email{liam@math.ubc.ca}
\thanks{AK is supported by an AMS-Simons travel grant. LW is supported by an NSERC discovery/accelerator grant.}

\author{Claudius Zibrowius}
\address{Faculty of Mathematics \\ University of Regensburg}
\email{claudius.zibrowius@posteo.net}

\begin{abstract} 
	There is a one-to-one correspondence between strong inversions on knots in the three-sphere and a special class of four-ended tangles.  We compute the reduced Khovanov homology of such tangles for all strong inversions on knots with up to 9 crossings, and discuss these computations in the context of earlier work by the second author.  In particular, we provide a counterexample to \cite[Conjecture~29]{Watson2017} as well as a refinement of and additional evidence for \cite[Conjecture~28]{Watson2017}.
\end{abstract}
\maketitle

\noindent The Brieskorn spheres $\Sigma(2,q,2nq\mp1)$ may be obtained by Dehn surgery on a torus knot in the three-sphere, namely, these are the integer homology spheres $S^3_{\nicefrac{\pm1}{n}}(T_{2,q})$ where $T_{2,q}$ is the positive $(2,q)$ torus knot.  These homology spheres admit Seifert fibrations, with base orbifold $S^2(2,q,2nq\mp1)$. Denoting by $\boldsymbol{\Sigma}(A,b)$ the two-fold branched cover of $A$ with branch set $b$, each of these manifolds admits two descriptions as a two-fold branched cover:

\[S^3_{\nicefrac{\pm1}{n}}(T_{2,q})\cong \boldsymbol{\Sigma}(S^3,T^*_{q,2qn\mp1})\cong\boldsymbol{\Sigma}(S^3,\tau(\pm\tfrac{1}{n}))\]

\parpic[r]{
 \begin{minipage}{35mm}
 \centering
 \setlength{\captionmargin}{0pt}
 \captionsetup{type=figure}
 \includegraphics[scale=1]{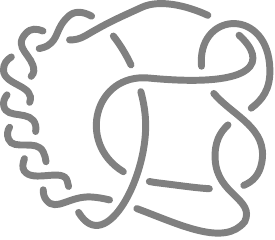}
 \captionof{figure}{The branch set $\tau(+1)$ when $q=5$}
 \label{fig:montesinos}
 \medskip
  \end{minipage}%
  }
\noindent This construction might be best termed as classical; for our purposes it is helpful to review the notation introduced in \cite{Watson2010}. The first of these two-fold branched covers results from an involution on the Seifert fibred space that preserves an orientation on the fibres---we refer to this as the Seifert involution. The second of these arises from the Montesinos involution, which reverses an orientation on the fibres; the branch set in question arises from the Montesinos trick, that is, by first constructing a tangle $(B^3,\tau)$ over which the exterior of $T_{2,q}$ is realized as a two-fold branched cover. We review the construction below as it is central to our enumeration of tangles. In particular, the branch set $\tau(\pm\tfrac{1}{n})$ is an explicit Montesinos link, which despite the notation depends on $q$. Generically, these two branch sets are distinct knots. For example, when $q=5$ and $n=1$ we have  
\[S^3_{+1}(T_{2,5})\cong \boldsymbol{\Sigma}(S^3,T^*_{5,9})\cong\boldsymbol{\Sigma}(S^3,\tau(+1))\]
and it can be calculated that 
\[\dim\Khr(T^*_{5,9})=57>15=\dim\Khr(\tau(+1))\]
where the reduced Khovanov homology is taken over the two-element field $\F$. (The Montesinos branch set is shown in Figure \ref{fig:montesinos}, and the construction of the tangle in this case is reviewed in Figure \ref{fig:cinqfoil-branch}.) This example serves to answer a question due to Ozsv\'ath in the negative: the total dimension of the mod 2 reduced Khovanov homology is not an invariant of two-fold branched covers \cite{Watson2010}. 

This paper focusses on tangles admitting knot exteriors as two-fold branched covers, and the immersed curves that arise as the reduced Khovanov invariants of these tangles. Coefficients are restricted to $\F$ throughout.  

\section{Strong inversions}  \label{sec:strong_inversions}
A knot $K$ in $S^3$ is invertible if it admits an isotopy exchanging a choice of orientation on $K$ with the reverse of this choice. A strong inversion is an inversion realized by an involution of the three-sphere. Note that invertible knots which are not strongly invertible exist \cite{Hartley}, but that when restricting to hyperbolic knots the two symmetries are equivalent. Following Sakuma \cite{Sakuma}, given a knot $K$ and a strong inversion $h$, we will call the pair $(K,h)$ a strongly invertible knot. Strongly invertible knots $(K,h)$ and $(K',h')$ are equivalent if there exists an orientation preserving homeomorphism $f$ on $S^3$ for which $f(K)=K'$ (so that $K$ and $K'$ are equivalent knots) and $h = f^{-1}\circ h'\circ f$. Alternatively, a strong inversion on a given knot may be viewed as the conjugacy class of an order-two element of the mapping class group of the exterior of the knot; this mapping class group is known as the symmetry group of the knot~$K$. For hyperbolic knots, this group is a subgroup of a dihedral group; see  \cite[Proposition 3.1]{Sakuma}. 

\labellist \tiny
\pinlabel $\mu$ at 5 23
\pinlabel $h$ at 57 72
\endlabellist
\parpic[r]{
 \begin{minipage}{35mm}
 \centering
 \medskip
 \setlength{\captionmargin}{0pt}
 \captionsetup{type=figure}
 \includegraphics[scale=1]{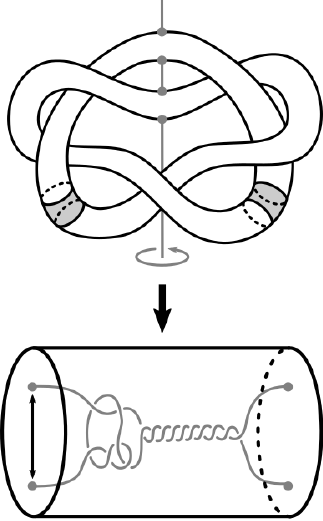}
 \captionof{figure}{The cinqfoil exterior as a two-fold branched cover}
 \label{fig:cinqfoil-branch}
  \end{minipage}%
  }
A strong inversion gives rise to an involution on a knot's exterior with one-dimensional fixed point set  meeting the boundary in exactly 4 points. Note that, according to the Smith conjecture, given an involution on the three-sphere with one-dimensional fixed point set, the fixed point set is unknotted; as a result the fixed point set of the strong inversion on the exterior is a pair of unknotted arcs. Taking the quotient of the order two action on the knot exterior gives rise to a three ball $B^3=(S^3\smallsetminus\nu(K)) / h$, and the image of the fixed point set of the strong inversion gives rise to a pair of properly embedded arcs $\tau=\operatorname{Im}(\operatorname{Fix}(h))$. Therefore, given a strongly invertible knot, the knot exterior may be viewed as a two-fold branched cover over a four-ended tangle $\tau$ in a three ball: 
$h \curvearrowright(S^3,K)$ so that $S^3\smallsetminus\nu(K) \cong \boldsymbol{\Sigma}(B^3,\tau)$. We refer to the tangle $T=(B^3,\tau)$ as the associated quotient tangle to a given knot with strong inversion; see Figure \ref{fig:cinqfoil-branch} for an explicit example.

\begin{figure}[b]
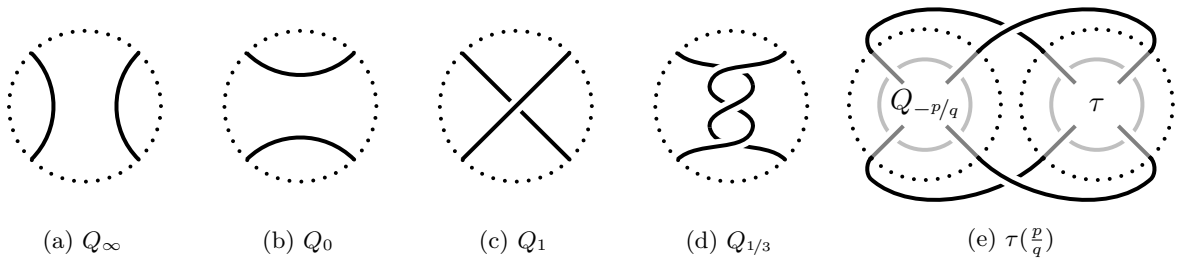

	\centering
	\begin{subfigure}{0.17\textwidth}
		\centering
		\(\ratixo\)
		\caption{\(Q_{\infty}\)}
		\label{fig:rat:ixo}
	\end{subfigure}
	\begin{subfigure}{0.17\textwidth}
		\centering
		\(\ratoxi\)
		\caption{\(Q_{0}\)}
		\label{fig:rat:oxi}
	\end{subfigure}
	\begin{subfigure}{0.17\textwidth}
		\centering
		\(\ratixi\)
		\caption{\(Q_{1}\)}
		\label{fig:rat:ixi}
	\end{subfigure}
	\begin{subfigure}{0.17\textwidth}
		\centering
		\(\ratixiii\)
		\caption{\(Q_{\nicefrac{1}{3}}\)}
		\label{fig:rat:ixiii}
	\end{subfigure}
	\begin{subfigure}{0.28\textwidth}
		\centering
		\(\tanglepairing\)
		\caption{\(\tau(\tfrac{p}{q})\)}
		\label{fig:rat:pairing}
	\end{subfigure}
	\caption{Examples of framed rational tangles~(a-d) and the rational closure of a four-ended tangle~\(\tau\)~(e)}
	\label{fig:rat}
\end{figure}

In this context, the natural notion of equivalence on tangles is homeomorphism of the pair $(B^3,\tau)$ where the boundary is fixed only set-wise. Note that this is, of course, more flexible than the requirement that the boundary be fixed point-wise as is perhaps more common when considering tangle diagrams. In order to be clear about the distinction, we will refer to this latter as a framed tangle. We remark that \cite[Definition 3]{Watson2010} gives a third notion of equivalence by introducing sutured tangles. This object will not play an explicit role here, however, consulting Figure \ref{fig:cinqfoil-branch} the reader may find that sutures are a helpful tool for tracking the image of the knot meridian $\mu$ in the quotient. 

From this point of view, there is a distinguished {\it trivial} tangle obtained as the quotient of a solid torus, that is, the associated quotient tangle for the trivial knot. The equivalence class of this tangle is the rational tangle. We will have need for choices of representatives \(Q_{\nicefrac{p}{q}}\) as described in Figure \ref{fig:rat}: For any \(\tfrac{p}{q}\in\QPI\), there is a framed tangle diagram \(Q_{\nicefrac{p}{q}}\). This choice will allow us to make use of the Montesinos trick: 
\[
S^3_{\nicefrac{p}{q}} (K) 
\cong 
\boldsymbol{\Sigma}(S^3,\tau(\tfrac{p}{q}))
\]
where \(\tfrac{p}{q}\)-surgery along the knot \(K\) corresponds to the knot \(\tau(\tfrac{p}{q})\) obtained by gluing the \(-\tfrac{p}{q}\)-rational tangle \(Q_{-\nicefrac{p}{q}}\) to the tangle \(\tau\) as in Figure~\ref{fig:rat:pairing}.

In particular, we are fixing a preferred representative for our associated quotient tangle once and for all:  Given $(K,h)$, this is the tangle \(T=(B^3,\tau)\) such that the rational closure $\tau(0)$ corresponds to surgery along the Seifert longitude and \(\tau(\infty)\) corresponds to \(S^3_{\infty} (K)\cong S^3\). The latter implies that \(\tau(\infty)\) is the unknot, as observed above. In fact, there is a bijection between equivalence classes of non-trivial strongly invertible knots $(K,h)$ and tangles $(B^3,\tau)$ for which $\tau(\infty)$ is unknotted (this follows from \cite[Theorem 2]{GL1989}; compare \cite[Proposition 9]{Watson2017}). As such, associated quotient tangles provide an invariant of strong inversions. 

Strong inversions give a means of enumerating interesting tangles, and we use this strategy below. Observe that, given a hyperbolic knot admitting a pair of distinct strong inversions $h_1$ and $h_2$, by Thurston's hyperbolic Dehn filling theorem, $S^3_{\nicefrac{p}{q}}(K)$ is hyperbolic for all but finitely many slopes \(\tfrac{p}{q}\). Moreover
\[
\boldsymbol{\Sigma}(S^3,\tau_1(\textstyle\frac{p}{q}))
\cong
S^3_{\nicefrac{p}{q}}(K)
\cong
\boldsymbol{\Sigma}(S^3,\tau_2(\textstyle\frac{p}{q}))
\]
and, generically, the branch sets $\tau_1(\tfrac{p}{q})$ and $\tau_2(\tfrac{p}{q})$ are distinct knots. It is a striking fact that for all $K$ with fewer than 9 crossings
\[\dim\Khr(\tau_1\textstyle(\frac{p}{q}))=\dim\Khr(\tau_2(\textstyle\frac{p}{q}))\]
that is, finding a negative answer to Ozsv\'ath's question in the hyperbolic setting is surprisingly difficult. This observation follows from Theorem \ref{thm:computation} below.

\begin{figure}[bt]
	\labellist \tiny
	\pinlabel $\mu$ at 85 27 \pinlabel $\mu$ at 218 27
	\pinlabel $\lambda+n_1\mu$ at 38 -4 \pinlabel $\lambda+n_2\mu$ at 167 -4
	\endlabellist
	\includegraphics{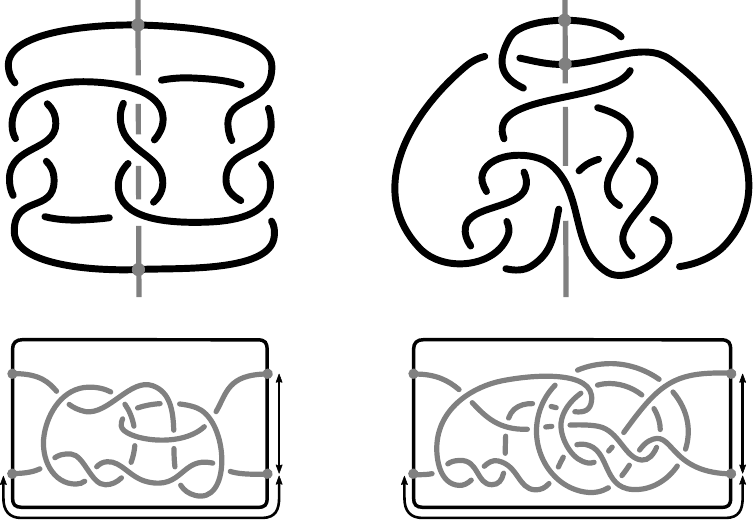}
	\caption{Two strong inversions on the knot $9_{46}$ (top row) together with the associate quotient tangle for each (bottom row). In order to present minimal crossing diagrams for the associated quotient tangles, this is the one place where we present something other than the preferred tangle representative: The first symmetry gives rise to a tangle $T_1$ shown with framing \(n_1=2\), while the second symmetry gives rise to a tangle $T_2$ with framing \(n_2=-6\).}\label{fig:9_46}
\end{figure}
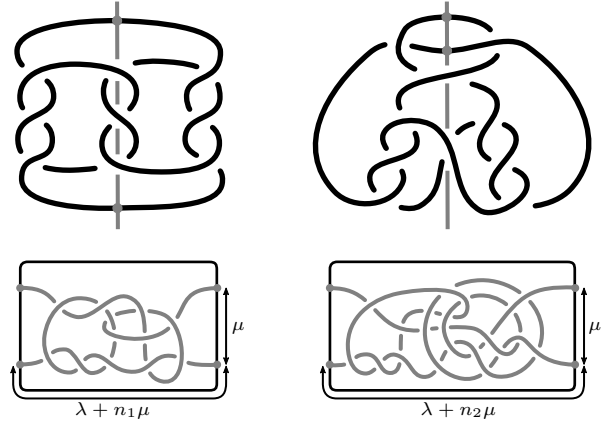

A range of examples of strong inversions and associated quotient tangles are considered in \cite{Watson2017}. We expand this list in a systematic way below, and highlight in particular the knot $9_{46}$ in the Rolfsen table; see Figure \ref{fig:9_46}. It is worth noting that these tangles distinguish the strong inversions in question: Ignoring one strand, the second tangle contains $8_{19}$ as a subknot, while the only subknots in the first tangle are trefoils. 

In the context of the discussion above, it seems natural to try and articulate a relative version of Ozsv\'ath's question. In particular, are there Khovanov-type tangle invariants that are invariants of the knot $K$ independent on the chosen strong inversion $h$?  We will see that $9_{46}$ dashes any hope of this (see Counterexample~\ref{counter}) and, furthermore, provides hyperbolic counterexamples to Ozsváth's original question (see Remark \ref{rmk:9_46}).

Interestingly, $9_{46}$ also makes a star appearance in the  recent work of Boyle and Issa \cite{BI}. In particular, this knot is slice, but it has non-zero equivariant four-genus; see \cite[Figure 14]{BI} and the surrounding discussion.

\section{\texorpdfstring{Review of the tangle invariants \(\Khr\) and \(\BNr\)}{Review of the tangle invariants Khr and BNr}}\label{sec:reviewKhr}

Let $T$ be a four-ended tangle in the three-ball \(B^3\).
In earlier work we interpreted Bar-Natan's tangle invariant $\KhTl{T}$ in terms of isotopy classes of immersed curves in the four-punctured sphere \(\partial B^3\smallsetminus \partial T\) \cite{KWZ,BarNatanKhT}. 
After choosing a distinguished tangle end $\ast$ of \(T\), the construction outputs two curve-valued invariants $\Khr(T)$ and $\BNr(T)$:
$$
T \quad 
\xmapsto{\text{\cite{BarNatanKhT}}}  
\quad 
\KhTl{T}   
\quad 
\xmapsto{\text{\cite{KWZ}}} 
\quad 
\Khr(T),~\BNr(T) \looparrowright \FourPuncturedSphereKh=\partial B^3\smallsetminus \partial T
$$
Strictly speaking, these immersed curves come equipped with local systems. However, in this paper, we will  suppress this subtlety, since over $\F$ the local systems are trivial in all known examples.

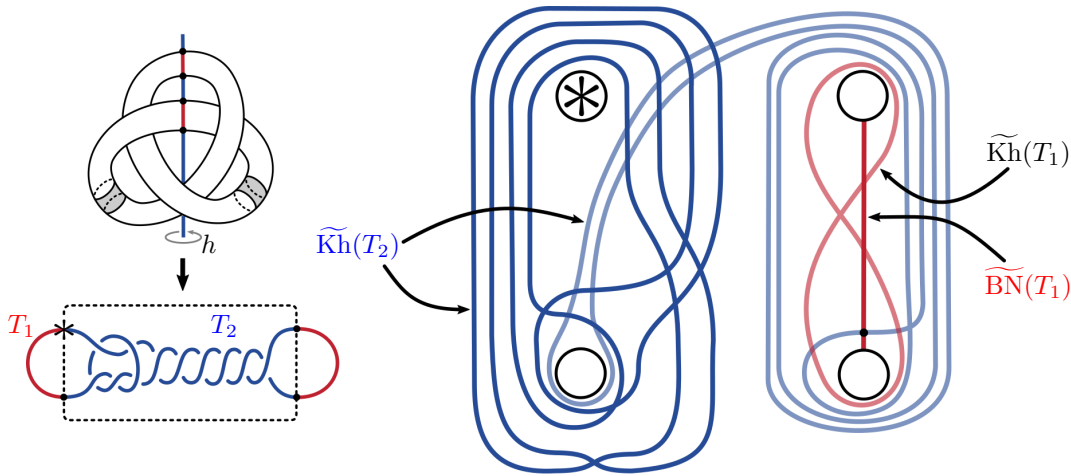
\begin{figure}[t]
\centering
\labellist 
\small
\pinlabel $h$ at 76.5 97
\pinlabel ${\color{blue}T_2}$ at 83 62
\pinlabel ${\color{red}T_1}$ at -2 62
\pinlabel ${\color{blue}\Khr(T_2)}$ at 139 97
\pinlabel ${\color{red}\BNr(T_1)}$ at 420 80
\pinlabel ${\Khr(T_1)}$ at 420 136
\endlabellist
\includegraphics[scale=0.9]{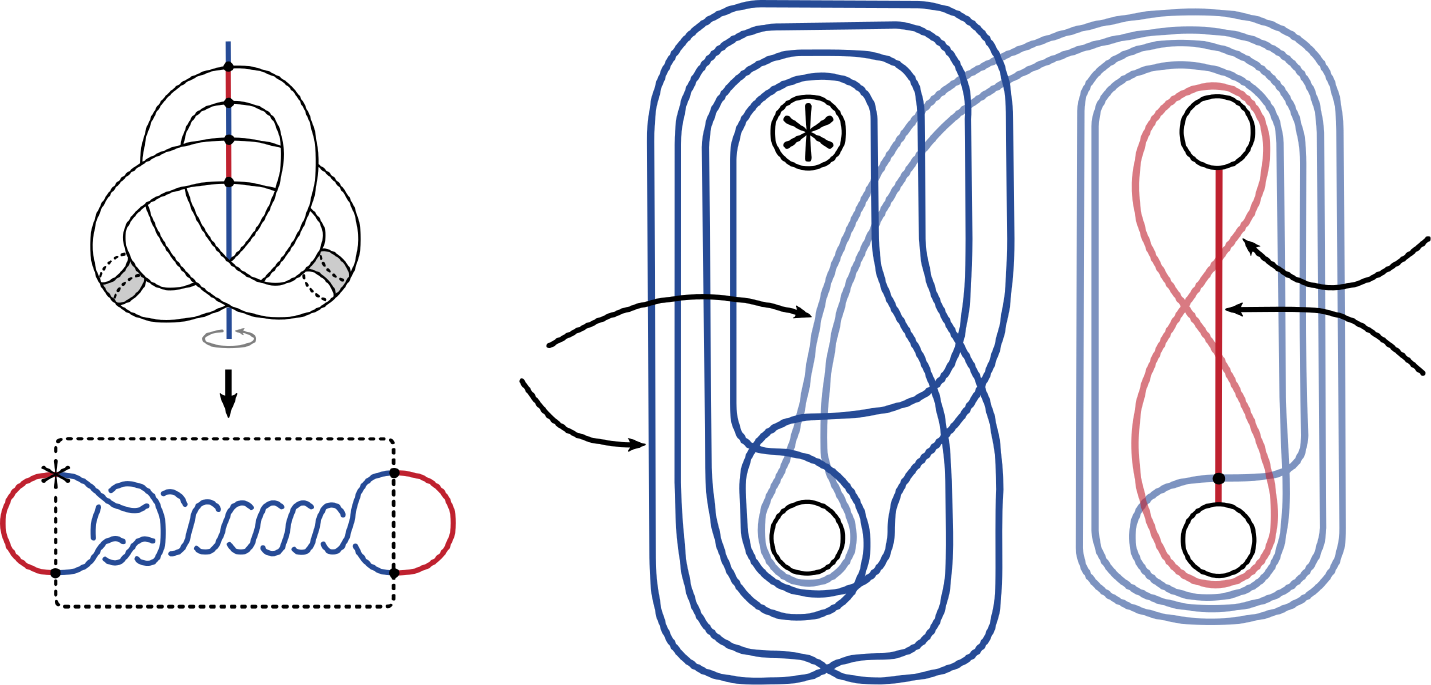}
\caption{Illustrated on the left is the trefoil with strong inversion, highlighting the decomposition of the unknot into tangles $\textcolor{red}{T_1}$ and $\textcolor{blue}{T_2}$, where $\textcolor{blue}{T_2}$ is the preferred representative of the associated quotient tangle to the strong inversion.  The corresponding intersection picture is shown on the right, where all of the relevant information has been projected to the front face of the sphere: The Lagrangian Floer homology $\HF(\textcolor{red}{\BNr(T_1)},\textcolor{blue}{\Khr(T_2)})=\F$ recovers the one-dimensional reduced Khovanov homology of the unknot.}\label{fig:pairing_example}
\end{figure}

The curve invariants enjoy the following gluing property \cite[Theorem~1.9]{KWZ}: Given a decomposition of a knot $K=T_1\cup T_2$ into four-ended tangles, 
the reduced Khovanov homology of the knot is recovered via Lagrangian Floer homology:
\begin{equation} \label{eq:pairing}
\Khr(K) ~ \cong ~ \HF(\mirror(\BNr(T_1)), \Khr(T_2))
\end{equation}
where $\mirror$
is the map identifying the two four-punctured spheres. 
An example illustrating this gluing formula is given in Figure~\ref{fig:pairing_example}.


Let us recall some basic facts about $\Khr(T)$ and $\BNr(T)$ from \cite[Section~6]{KWZ}: Both invariants may consist of multiple components. A component is either an immersion of a circle or an immersion of an interval. In the first case we call a component compact; in the second, non-compact. The invariant $\BNr(T)$ contains exactly one non-compact component (unless the tangle \(T\) contains some closed component, which is not the case in the present context), whereas $\Khr(T)$ consists only of compact components. 

These curves become easier to manage when considered in a certain covering space of \(\FourPuncturedSphereKh\), namely the planar cover that factors through the toroidal two-fold cover: 
$$(\R^2 \smallsetminus \Z^2) \to (T^2 \smallsetminus 4\pt) \to \FourPuncturedSphereKh$$
Given an immersed curve $c \looparrowright \FourPuncturedSphereKh$, denote by $\tc \looparrowright \R^2 \smallsetminus \Z^2$ an arbitrary lift of $\c$. 
A curve \(c\) is called \emph{linear} if there exists some \(\tfrac{p}{q}\in\QPI\) such that for every open neighbourhood \(U\) of \(\tfrac{p}{q}\) in \(\QPI\), there is a curve $\c_U$ homotopic to \(c\) with the property that the slopes \(\tc_U'(t)\) of the lift $\tc_U$ are contained in \(U\).
If there exists such a \(\tfrac{p}{q}\in\QPI\), it is unique, and we call it the \emph{slope} of \(c\). 
We call a linear curve \(\c\) \emph{rational} if there exists some neighbourhood of the special puncture \(\ast\) of \(\FourPuncturedSphereKh\) that $\c_U$ can be chosen to avoid for all \(U\) simultaneously. 
Otherwise, we call \(\c\) \emph{special}.


\begin{figure}[t]
\vspace{0.5cm}
\centering
\labellist 
\pinlabel $\overbrace{\hspace{3.2cm}}^{\text{2n punctures}}$ at 158 42.5
\pinlabel $\overbrace{\hspace{2.2cm}}^{\text{n punctures}}$ at 172 13.5
\small
\pinlabel $\ts_n(0)$ at 20 45
\pinlabel $\tr_n(0)$ at 90 15.5
\endlabellist
\includegraphics[scale=1]{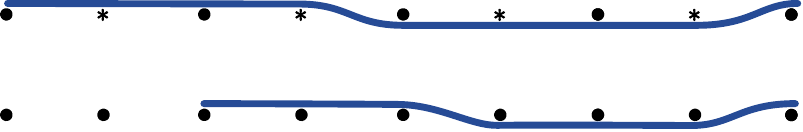}
\caption{Curves $\{\ts_n(0)\}_{n\geq 1}$ and $\{\tr_n(0)\}_{n\geq 1}$}\label{fig:cover}
\end{figure}

Figure~\ref{fig:cover} shows some examples of lifts of linear curves of slope 0. The curves $\ts_n(0)$ define a family  of special curves \(\s_n(0)\), \(n\in\Z^{>0}\); the curves \(\tr_n(0)\) define a family  of rational curves  \(\r_n(0)\), \(n\in\Z^{>0}\).  For every \(\tfrac{p}{q}\in\QPI\), we define the linear curves $\s_n(\frac p q)$ and $\r_n(\frac p q)$ of slope \(\tfrac{p}{q}\) as the images of $\s_n(0)$ and $\r_n(0)$, respectively, under the action of
\[
\begin{bmatrix*}[c]
	q & r \\
	p & s
\end{bmatrix*}
\]
considered as an element the mapping class group $\operatorname{Mod}(\FourPuncturedSphereKh) \cong PSL(2,\Z)$, where \(qs-pr=1\). 

\begin{example}
	The curve $\Khr(T_2)$ from Figure~\ref{fig:pairing_example} consists of the special component $\s_2(\infty)$ and the rational component $\r_1(4)$. 
	Furthermore, $\BNr(T_1)=\{\a(\infty)\}$ is the red vertical arc, while $\Khr(T_1)=\{\r_1(\infty)\}$ consists of the figure-eight curve which lies in a small neighbourhood of this arc. More generally, justifying the terminology, naturality of the invariants under the mapping class group action implies that $\BNr(Q_{\nicefrac{p}{q}})=\{\a(\frac p q)\}$ and $\Khr(Q_{\nicefrac{p}{q}})=\{\r_1(\frac p q)\}$---so, rational tangles have rational invariants.
\end{example}

Computations suggest that the geography of components of $\Khr(T)$ is highly contrained:

\begin{enumerate}
\item Every component of \(\Khr(T)\) is linear.

\item Every special component of \(\Khr(T)\) is equal to \(\s_n(\frac p q)\) for some \(n \geq 1\) and \(\frac p q \in \QPI\).

\item Every rational component of \(\Khr(T)\) is equal to \(\r_n(\frac p q)\) for some \(\frac p q \in \QPI\).
\end{enumerate}
We expect to establish these restrictions in the near future~\cite{KWZ_thin}; in the context of Heegaard Floer homology, analogous properties are known for the tangle invariant \(\HFT(T)\)~\cite{pqSym}.
Note, however, that components of the invariant $\BNr(T)$ are known to be much more complicated even when restricting to only compact components.

We conclude this section with a computation of the Lagrangian Floer homology of a simple arc with the special curves \(\s_n(\infty)\), which will play a central role in the next section. 

\begin{definition}
	Define \(\sKappa\) as the four-dimensional bigraded vector space supported in \(\delta\)-grading 0 and quantum gradings \(-5\), \(-1\), \(+1\), and \(+5\). 
\end{definition}

\begin{figure}
	\[
	\begin{tikzcd}[blue]
	\GGzqh{\DotVblue}{+\frac{1}{2}}{-5}{-3}
	\arrow{d}{D}
	\arrow{r}{S}
	&
	\GGzqh{\DotHblue}{0}{-4}{-2}
	\arrow{r}{D}
	&
	\GGzqh{\DotHblue}{0}{-2}{-1}
	\arrow{r}{S^2}
	&
	\GGzqh{\DotHblue}{0}{0}{0}
	\arrow{r}{D}
	&
	\GGzqh{\DotHblue}{0}{2}{1}
	\arrow{r}{S}
	&
	\GGzqh{\DotVblue}{-\frac{1}{2}}{3}{2}
	\arrow{d}{D}
	\\
	\GGzqh{\DotVblue}{+\frac{1}{2}}{-3}{-2}
	\arrow{r}{S}
	&
	\GGzqh{\DotHblue}{0}{-2}{-1}
	\arrow{r}{D}
	&
	\GGzqh{\DotHblue}{0}{0}{0}
	\arrow{r}{S^2}
	&
	\GGzqh{\DotHblue}{0}{2}{1}
	\arrow{r}{D}
	&
	\GGzqh{\DotHblue}{0}{4}{2}
	\arrow{r}{S}
	&
	\GGzqh{\DotVblue}{-\frac{1}{2}}{5}{3}
	\end{tikzcd}
	\]
	\caption{The symmetrically bigraded type D structure corresponding to the special curve \(\textcolor{blue}{\s_{2n}(\infty)}\) for \(n=1\). The type~D structures for \(n>1\) look similar. The total number of generators in idempotent \(\protect\DotHblue\) is equal to \(8n\). The absolute bigrading is fixed by requiring that the minimal and maximal homological and quantum gradings are \(\pm(2n+1)\) and \(\pm(4n+1)\), respectively.}
	\label{fig:complex_for_special}
\end{figure}
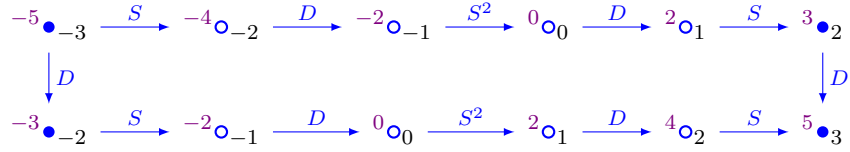

\begin{figure}
	\labellist \tiny
	\pinlabel \QGrad{$0$}  at 129 46 \pinlabel $0$ at 135.5 37.5
	\pinlabel \QGrad{$5$}  at 129 69 \pinlabel $3$ at 135.5 61
	\pinlabel \QGrad{$3$}  at 129 82 \pinlabel $2$ at 135.5 73.5
	\pinlabel \QGrad{-$3$} at 33 72 \pinlabel -$2$ at 41 63
	\pinlabel \QGrad{-$5$} at 33 30 \pinlabel -$3$ at 41 21
	\pinlabel \QGrad{-$2$} at 107 182 \pinlabel -$1$ at 117 175
	\pinlabel \QGrad{-$4$} at 128 182 \pinlabel -$2$ at 138 175
	\pinlabel \QGrad{$0$}  at 150 182 \pinlabel $0$ at 158 175
	\pinlabel \QGrad{$2$}  at 171 182 \pinlabel $1$ at 178 175
	\pinlabel \QGrad{$0$}  at 108.5 147 \pinlabel $0$ at 117 139.5
	\pinlabel \QGrad{-$2$} at 128 147 \pinlabel -$1$ at 139 139.5
	\pinlabel \QGrad{$2$}  at 150 147 \pinlabel $1$ at 158 139.5
	\pinlabel \QGrad{$4$}  at 171 147 \pinlabel $2$ at 179 139.5
	\small
	\pinlabel $\QGrad{q^5}\delta^{-\frac{1}{2}}h^3$ at 295 116.5
	\pinlabel $\QGrad{q^{1}}\delta^{-\frac{1}{2}}h^{1}$ at 408 116.5
	\pinlabel $\QGrad{q^{-5}}\delta^{-\frac{1}{2}}h^{-2}$ at 295 -3
	\pinlabel $\QGrad{q^{-1}}\delta^{-\frac{1}{2}}h^{0}$ at 408 -3
	\tiny
	\pinlabel $S$ at 85 145
	\pinlabel $D$ at 61 145
	\pinlabel $S$ at 61 172
	\pinlabel $D$ at 87 172
	\endlabellist
	\includegraphics[scale=0.9]{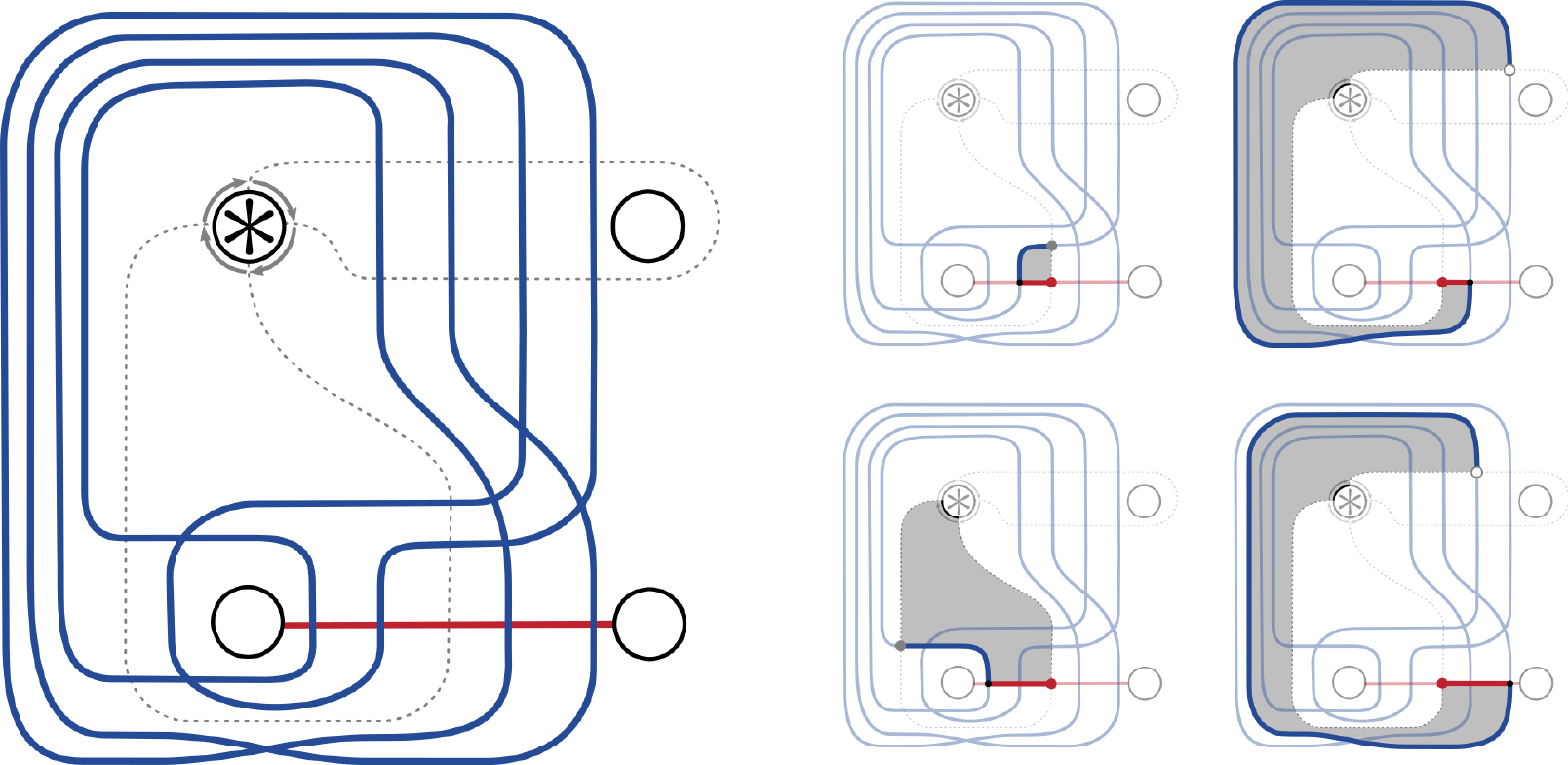}\medskip
	\caption{Computing the gradings associated with the four-dimensional vector space \(\HF(\textcolor{red}{\a(0)},\textcolor{blue}{\s_2(\infty)})\), according to Lemma \ref{lem:kappa_from_s_n}}\label{fig:pairing_special_with_arc}
\end{figure}

\begin{lemma}\label{lem:kappa_from_s_n}
	Let \(\textcolor{red}{\a(0)}=\BNr(\No)\) and \(\textcolor{blue}{\s_2(\infty)}\) as in Figure~\ref{fig:pairing_special_with_arc}. 
	Then 
	\[
	\HF(\textcolor{red}{\a(0)},\textcolor{blue}{\s_2(\infty)})
	\cong
	q^{0}
	\delta^{-\tfrac{1}{2}}
	h^{\tfrac{1}{2}}
	\sKappa
	\]
	More generally, for any positive integer \(n\), 
	\[
	\HF(\textcolor{red}{\a(0)},\textcolor{blue}{\s_{2n}(\infty)})
	\cong
	\bigoplus_{i=1}^n
	q^{4(2i-n-1)}
	\delta^{-\tfrac{1}{2}}
	h^{\tfrac{1}{2}+2(2i-n-1)}
	\sKappa
	\]
	In both cases, the isomorphism holds as absolutely bigraded vector spaces if  \(\textcolor{blue}{\s_{2n}(\infty)}\) is symmetrically bigraded in the sense of Figure~\ref{fig:complex_for_special}. 
\end{lemma}

\begin{proof}
	The total dimension of \(\HF(\textcolor{red}{\a(0)},\textcolor{blue}{\s_{2n}(\infty)})\) is \(4n\), since the minimal number of intersection points between the two curves is \(4n\).
	To compute the bigrading on the Lagrangian Floer homology, we need to recall some facts from \cite[Section~7]{KWZ}. We will focus on the case \(n=1\), the computation for \(n>1\) is similar.
	
	The invariants $\BNr(T)$ and $\Khr(T)$ are topological interpretations of algebraic invariants \(\DD(T)\) and \(\DD_1(T)\), respectively, which are type~D structures over the algebra \(\BNAlgH\). For example, the curve \(\textcolor{blue}{\s_2(\infty)}\), which is the blue curve on the left of Figure~\ref{fig:pairing_special_with_arc}, corresponds to the type~D structure in Figure~\ref{fig:complex_for_special}; the arc \(\textcolor{red}{\a(0)}\) corresponds to the one-generator complex \([\GGzqh{\DotVred}{0}{0}{0}]\). 
	To translate between these two viewpoints, we need to fix a parametrization of the four-punctured sphere \(\FourPuncturedSphereKh\), which is indicated by the two grey dashed arcs in Figure~\ref{fig:pairing_special_with_arc}.   Briefly,  intersection points of these arcs with a given curve give rise to the generators of the corresponding type~D structure and paths between intersection points determine the differential; for more details, see \cite[Example~1.6]{KWZ}. The generators of the type~D structures carry a bigrading, that is a quantum grading \(q\) and a homological grading \(h\). We indicate these gradings on generators \(\bullet\) and the corresponding intersection points by super- and subscripts like so: \(\GGzqh{\bullet}{}{q}{h}\). Then the bigrading of an intersection point \(\bullet\) generating \(\HF(\textcolor{red}{\a(0)},\textcolor{blue}{\s_{2}(\infty)})\) is computed by considering a path on \(\textcolor{red}{\a(0)}\cup\textcolor{blue}{\s_{2}(\infty)}\) which starts at \(x=\GGzqh{\DotVred}{0}{0}{0}\), turns right at the intersection point \(\bullet\) and ends at the nearest intersection point \(y\) of  \(\textcolor{blue}{\s_{2}(\infty)}\) with the parametrization. 
	These paths are illustrated on the right of Figure~\ref{fig:pairing_special_with_arc}. 
	Each of these paths in \(\FourPuncturedSphereKh\) is homotopic relative to the parametrizing arcs to some path \(\gamma\) on the boundary of the special puncture; these homotopies are indicated on the right of Figure~\ref{fig:pairing_special_with_arc} by the shaded disks.  We set \(\QGrad{q}(\gamma)=\QGrad{0},\QGrad{-1},\QGrad{-2}\), depending on whether the path \(\gamma\) is constant, equal to \(S\), or equal to \(D\), respectively. (More generally, these paths correspond to some algebra elements in \(\BNAlgH\), whose quantum gradings define \(\QGrad{q}(\gamma)\).) Then
	\[
		h(\bullet)=h(y)-h(x)
		\quad
		\text{and}
		\quad 
		\QGrad{q}(\bullet)=\QGrad{q}(y)-\QGrad{q}(x)+\QGrad{q}(\gamma)
	\]
	Finally, we compute the \(\delta\)-grading from the identity \(\delta+h=\frac 1 2 \QGrad{q}\). 
\end{proof}

\section{The Khovanov homology of strong inversions}

A strongly invertible knot $(K,h)$ gives rise to a four-ended associated quotient tangle $T=(B^3,\tau)$.  (When we have need for it, we will use the notation $T_{K,h}$ to remember the dependence on the strongly invertible knot.) We use the same preferred framing as in Section~\ref{sec:strong_inversions}; in particular,  $\tau(\infty)$ is the unknot. This imposes strong restrictions on the curve invariant $\Khr(T)$, namely, it implies that
$$
\Khr(\tau(\infty))
=
\HF(\BNr(\Ni),\Khr(T))
=
\HF(\a(\infty),\Khr(T))
=
\F
$$
In other words, \(\Khr(T)\) is homotopic to a multicurve which intersects the arc \(\a(\infty)\) only once, as is the case for the example $\Khr(T_2)$ from Figure~\ref{fig:pairing_example}. 
In the light of the geography restrictions, this means that $\Khr(T)$ should be a union of special components of slope $\infty$ together with a single rational component of an integer slope. 
In fact, with the help of a computer \cite{khtpp}, we show the following: 

\begin{theorem}\label{thm:computation}
	Let \(h\) be a strong inversion on a knot \(K\subset S^3\) with at most 9 crossings. Then
	$$
	\Khr(T_{K,h})
	=
	\r_1(k)
	\cup
	\underbrace{%
		\s_2(\infty)
		\cup
		\dots
		\cup
		\s_2(\infty)
	}_n
	$$ 
	for some integers \(k\in4\Z\) and \(n\in\Z^{>0}\). \qed
\end{theorem}

\begin{table}[t]
	\begin{multicols}{5}
		\input{computations/summary}
	\end{multicols}
	\caption{
	The ungraded reduced Khovanov homology of tangles \(T_{K,h}\) associated with strong inversions \(h\) on knots \(K\) with up to 9 crossings. The first columns specify the strong inversions, the second columns specify the slope \(k\) of the single rational component \(\r_1(k)\) in \(\Khr(T_{K,h})\), and the third columns give the number \(n\) of special components, all of which are equal to \(\s_2(\infty)\). 
	}
	\label{tab:summary}
\end{table}

	Table~\ref{tab:summary} shows the ungraded invariants for all pairs \((K,h)\) from Theorem~\ref{thm:computation}; the bigraded invariants are listed in Section~\ref{sec:kappa_table}. Note that, in order to match Sakuma's table~\cite{Sakuma}, the first row in Table~\ref{tab:summary} describes the left-hand trefoil and its rational component has slope $-4$, while Figure~\ref{fig:pairing_example} depicts the right-hand trefoil and its rational component has slope $4$. Among tabulated knots through $9$ crossings, there are $57$ knots that admit two distinct strong inversions. For all but one of these knots with two distinct strong inversions, the ungraded invariants \(\Khr(T)\) agree; these pairs of strong inversions are indicated in Table~\ref{tab:summary} by the superscript \(\star\).

\begin{remark}\label{rmk:9_46}
	The knot \(9_{46}\) is the only knot in this collection whose strong inversions \(9_{46}^1\) and \(9_{46}^2\) can be distinguished by their ungraded invariants \(\Khr(T)\), as shown in the two highlighted rows in Table~\ref{tab:summary}. Here, we follow the same numbering convention for strong inversions as in~\cite{Sakuma}. Note that both the slope \(k\) and the number \(n\) of special components \(\s_2(\infty)\) distinguish \(9_{46}^1\) and \(9_{46}^2\). With this example in hand, one can check that $\dim\Khr(\tau_1(n))\ne \dim\Khr(\tau_2(n))$ by applying the pairing formula~\eqref{eq:pairing}, despite the fact that $\mathbf{\Sigma}(S^3,\tau_1(n))\cong\mathbf{\Sigma}(S^3,\tau_2(n))$ by construction. Calculations for $\dim\Khr(\tau_i(\frac{p}{q}))$ with $i=1,2$ can be obtained as well with a little more patience. 
	\end{remark}

As stated in Section~\ref{sec:reviewKhr}, we expect \(\Khr(T)\) to contain only linear curves; it would be interesting to give a geometric interpretation for the slopes of these curves. For quotient tangles \(T_{K,h}\) of strong inversions \((K,h)\), the existence of the unknot closure implies that the slope of any special curve is fixed and the slope of the rational component is an integer. As we can see from the second column in Table~\ref{tab:summary}, this integer may be non-zero; in other words, the slope of the rational component of \(\Kh(T_{K,h})\) need not agree with the slope given by the rational longitude of the knot~\(K\). Our computations raise the following questions: 

\begin{questions}\label{q:rationals}
	Is there a geometric/topological meaning of the slope of the rational component of \(\Khr(T_{K,h})\)? 	Is the slope always divisible by 4?
\end{questions}

We now focus on the special components:

\begin{definition}\label{def:kappa}
	Given a knot with a strong inversion \((K,h)\), let \(\s(K,h)\) be the set of special components of \(\Khr(T_{K,h})\). Let \(\a(0)=\BNr(\No)\) as in Figure~\ref{fig:pairing_special_with_arc}. 
	Then define 
	\[
	\Kappa(K,h)
	\coloneqq 
	\HF(\a(0),\s(K,h))
	\]
\end{definition}

One can show that \(\Kappa(K,h)\) agrees as a relatively bigraded invariant with an invariant of the same name that was defined by the second author as a certain finite-dimensional quotient of the inverse limit of \(\displaystyle\lim_{\longleftarrow}\Khr(\tau(n))\) associated with the maps \(\Khr(\tau(n+1))\rightarrow\Khr(\tau(n))\) induced by resolving a crossing \cite{Watson2017}. 

\begin{counterexample}\label{counter}
Returning to the pair of strong inversions on the knot $9_{46}$, we calculate that \[\dim\Kappa(K,h_1)=16<28=\dim\Kappa(K,h_2)\] contrary to \cite[Conjecture~29]{Watson2017}. This conjecture was originally posed in the hope that a relative variant of Ozsv\'ath's question might help to explain the examples described in \cite{Watson2010}, which depend heavily on the Seifert structure. Comparing with Remark \ref{rmk:9_46}, it is worth noting a simple  relationship between $\Kappa(K,h)$ and $\Khr(\tau(n))$ when restricting to the preferred framing: According to Definition \ref{def:kappa},
\[
\dim\Khr(\tau(n))
=
\dim\Kappa(K,h) + \dim \HF(\a(n),\r_1(k))
\]
where \(k\) is as in Theorem~\ref{thm:computation}. 
Note that the second summand is the dimension of the reduced Khovanov homology of the \((2,n-k)\)-torus link.
\end{counterexample}

On the other hand, certain structural properties appear to persist: 

\begin{conjecture}[{\cite[Conjecture~28]{Watson2017}}]\label{conj:geography_kappa_original}
	For any strong inversion \((K,h)\), the vector space \(\Kappa(K,h)\) is a direct sum of copies of \(\sKappa\). In particular, \(\dim\Kappa(K,h)\) is divisible by 4. 
\end{conjecture}

By Lemma~\ref{lem:kappa_from_s_n}, the vector space \(\sKappa\) is precisely the Lagrangian Floer homology of the arc \(\a(0)\) and the special curve \(\s_2(\infty)\).
So a strong inversion \((K,h)\) satisfies Conjecture~\ref{conj:geography_kappa_original} if all special components of \(\Khr(T_{K,h})\) are equal to \(\s_2(\infty)\) up to a grading shift. In particular, the calculations summarized in Table~\ref{tab:summary} verify Conjecture~\ref{conj:geography_kappa_original} for all strong inversions on tabulated knots through 9 crossings. 

We remark that Conjecture~\ref{conj:geography_kappa_original} implies that the slope of the rational component of \(\Khr(T_{K,h})\) should be divisible by 4 (Questions~\ref{q:rationals}, second part). This can be seen as follows. First note that  
$$|H_1(\mathbf{\Sigma}(S^3,\tau(0)))|=|H_1(S^3_{0}(K))|=\infty$$
hence \(\det(\tau(0))=0\). 
The determinant of a knot can also be computed by evaluating the Jones polynomial at \(-1\), which agrees with the Euler characteristic of Khovanov homology with respect to the $\delta$-grading:
$$0=\det(\tau(0))=|V(-1)|=|\chi_\delta \Khr(\tau(0))|=|\chi_\delta \HF(\mathbf{a}(0),\Khr(T))|$$
This means that the $\delta$-graded intersections of $\mathbf{a}(0)$ with special components cancel out intersections with the rational component. Conjecture~\ref{conj:geography_kappa_original} implies that intersections with special components come in groups of $4$ concentrated in a single $\delta$-grading, and so the number of intersections with  the rational component---also concentrated in a single $\delta$-grading---should be divisible by 4. 

Due to Lemma~\ref{lem:kappa_from_s_n}, the following is a refinement of Conjecture~\ref{conj:geography_kappa_original}:  

\begin{conjecture}\label{conj:even_no_specials}
	Given a strongly invertible knot \((K,h)\), any special component of \(\Khr(T_{K,h})\) is equal to \(\s_{2n}(\infty)\) for some \(n>0\). 
	In particular, all special components have even length.
\end{conjecture}

The special components in all examples that we have tabulated are equal to \(\s_2(\infty)\), but there exist strongly invertible knots \((K,h)\) for which \(\Khr(T_{K,h})\) contains special components \(\s_{2n}(\infty)\) with \(n\neq1\).

As we can see from Table~\ref{tab:summary}, the ungraded invariant \(\Kappa(K,h)\) cannot tell all strong inversions apart. However, the \emph{absolutely bigraded} invariant \(\Kappa(K,h)\) can distinguish any pair of strong inversion for knots up to 9 crossings. 



\begin{question}[{\cite[Question 19]{Watson2017}}]
	Is there any hyperbolic knot \(K\subset S^3\) with two distinct strong inversions that cannot be distinguished by the absolutely bigraded invariant \(\Kappa\)?
\end{question}

The restriction to hyperbolic is important here, and was omitted in error in the statement of \cite[Question 19]{Watson2017}. Note that the question is not interesting for torus knots, as these admit unique strong inversions. However, it is possible to generate satellite knots with strong inversions whose associated quotient tangles differ by mutation on a sub-tangle. Owing to the insensitivity of Khovanov homology to mutation, such symmetries cannot be separated. 

\begin{remark}
	We have focussed entirely on four-ended tangles admitting an unknot closure. 
	We expect Conjecture~\ref{conj:even_no_specials} to hold more generally: 	
	For any four-ended tangle \(T\), any special component of \(\Khr(T)\) should be equal to \(\s_{2n}(\tfrac{p}{q})\) for some \(\tfrac{p}{q}\in\QPI\) and \(n>0\). 
\end{remark}

\begin{remark}
In this paper, we have worked exclusively over the field \(\fieldTwoElements\) of two elements; all of the above questions should be read with this coefficient system in place. We conclude with a comment about other fields of coefficients. 
Unlike knot Floer homology, Khovanov homology is known to behave very differently over fields of different characteristic. Many of the conjectures above fail if we work over a field different from \(\fieldTwoElements\).
	For example, consider the second strong inversion \(7_4^2\) on the knot \(7_4\) in Sakuma's table \cite{Sakuma}.  Its invariant \(\Khr(T_{7_4^2};\Z/3)\) consists of a single rational component \(\r_1(6)\), two special components \(\s_2(\infty)\), but also five special components \(\s_1(\infty)\). 
	So the corresponding version of Conjecture~\ref{conj:geography_kappa_original} over \(\Z/3\) is false. Also note that the slope of the rational component is no longer divisible by 4. 
\end{remark}

\section{A table of invariants}
\label{sec:kappa_table}
In the following, we list the absolutely bigraded invariants \(\Kappa\) for all strong inversions whose underlying knots have at most 9 crossings. 
We do this as follows:
We subdivide the drawing plane into a grid and associate with each point on the plane a bigrading, namely the homological grading is equal to the \(x\)-coordinate and the \(\delta\)-grading is equal to the \(y\)-coordinate. 
Moving a point right and upwards increases both gradings. For each pair \((K,h)\), we can write
\[
\Kappa(K,h)
=
\bigoplus_{i}h^{a_i}\delta^{b_i}\sKappa
\]
since all special components of the tangle invariants have length 2. Then, for each \(i\), we place a grey box centred at the point \((a_i,b_i)\). If multiple boxes line up at the same bigrading, we add a label which indicates the number of such boxes. Finally, the absolute bigrading is specified by the two numbers in the bottom left corner, which indicate the absolute coordinate of the bottom-and-left-most intersection point of the grid: The top-left number is the \(\delta\)-grading, the bottom-right number is equal to the homological grading.   
For example, 
\[
\Kappa(4_{1}^{1})
=
h^{2.5}\delta^{5.5}\sKappa
\oplus
h^{6.5}\delta^{4.5}\sKappa
\] 
The superscripts \(\star=1,2\) indicate the strong inversions, using the same numbering convention as in \cite{Sakuma}. 
We follow the grading conventions in~\cite{KWZ}; when comparing these computations to \cite{Watson2017}, note that our \(\delta\)-grading has the opposite sign, and our quantum grading is twice as large. 

\begin{remark}
	There is a mistake in the diagram for \(9_{43}\) in \cite{Sakuma}, where the diagram given is an alternating knot of determinant 43 and Rasmussen's \(s\)-invariant \(\pm2\). The only such knots with $\leq 9$ crossings sharing these properties are \(9_{21}\) and \(9_{22}\), according to knotinfo \cite{knotinfo}. By comparing the \(\Kappa\)-invariant of the quotients, we see that Sakuma's diagram shows \(9_{22}\); The actual knot \(9_{43}\) can be obtained from this diagram by changing the two crossings in the centre of the diagram, close to the two points on the knot that stay fixed under the strong inversion. Its determinant is \(13\) and \(s=\pm4\). The only other knot of no more than $9$ crossings with these properties is \(7_3\), but its strong inversions have different \(\Kappa\) than the ones we compute for \(9_{43}\). 
\end{remark}

\begin{multicols}{2}
	\noindent
	\input{computations/gradings}
\end{multicols}

\begin{small}
	\pdfbookmark[1]{Acknowledgements}{Acknowledgements}
	\noindent\textbf{Acknowledgements.}
	The computer program \cite{khtpp} used in this work represents a major overhaul to an earlier piece of software, developed as part of a UBC undergraduate research project with Gurkeerat Chhina. 
\end{small}

\newcommand*{\arxiv}[1]{\href{http://arxiv.org/abs/#1}{ArXiv:\ #1}}
\newcommand*{\arxivPreprint}[1]{\href{http://arxiv.org/abs/#1}{ArXiv preprint #1}}

\bibliographystyle{alpha}
\bibliography{main}

\end{document}

%% file: notation_commands.tex
\renewcommand{\geq}{\geqslant}
\renewcommand{\leq}{\leqslant}

\newcommand{\fieldTwoElements}{\mathbb{F}}
\newcommand{\F}{{\mathbb{F}}}

\newcommand{\Z}{\mathbb{Z}}
\newcommand{\ZZ}{\mathbb{Z}/2}
 
\newcommand{\R}{\mathbb{R}}

\newcommand{\BNAlgH}{\mathcal{B}}

\newcommand{\pt}{\text{pt}}

\newcommand{\QGrad}[1]{{\textcolor{violet}{#1}}}
\newcommand{\DeltaGrad}[1]{{\textcolor{darkgreen}{#1}}}
\newcommand{\HomGrad}[1]{{\textcolor{black}{#1}}}
\newcommand{\AllGrad}[3]{\QGrad{q^{#1}}\DeltaGrad{\delta^{#2}}\HomGrad{h^{#3}}}
\newcommand{\GGdqh}[4]{\prescript{\QGrad{#3}}{}{#1}^\DeltaGrad{#2}_\HomGrad{#4}}

\newcommand{\GGzqh}[4]{\prescript{\QGrad{#3}}{}{#1}_\HomGrad{#4}}

\newcommand{\HF}{\operatorname{HF}}

\newcommand{\Kappa}{\varkappa}
\newcommand{\sKappa}{V_\varkappa}

\DeclareMathOperator{\HFT}{HFT}

\newcommand{\KhTl}[1]{[\![ #1 ]\!]_{/l}} 
\DeclareMathOperator{\Kh}{Kh}
\DeclareMathOperator{\BN}{BN}
\newcommand{\Khr}{\widetilde{\Kh}}
\newcommand{\BNr}{\widetilde{\BN}}

\newif\ifA
\newif\ifi
\newif\ifm
\newif\ifd
\newcommand{\thth}[1]{
	\operatorname{%
		\ifd%
			\partial%
		\fi%
		\ifA%
			\ifi%
				\mathring{\mathrm{A}}%
			\else%
				\mathrm{A}%
			\fi%
		\else%
			\ifi%
				\mathring{\Theta}%
			\else%
				\Theta%
			\fi%
		\fi%
		\ifm%
			\ifA%
				^{\!\mirror}%
			\else%
				^{\mirror}%
			\fi%
		\fi%
		_{#1}%
	}%
}
\newcommand{\HFsup}{\mathrm{HF}}
\newcommand{\Khsup}{\mathrm{Kh}}
\newcommand{\Gsup}{\mathit{G}}
\newcommand{\makemirror}[2]  {\newcommand{#1}{\mtrue#2\mfalse}}
\newcommand{\makeALink}[2]   {\newcommand{#1}{\Atrue#2\Afalse}}
\newcommand{\makeInterior}[2]{\newcommand{#1}{\itrue#2\ifalse}}
\newcommand{\makeBoundary}[2]{\newcommand{#1}{\dtrue#2\dfalse}}

\newcommand{\Thin}{\thth{}}
\newcommand{\ThinZ}{\thth{\Z}}
\newcommand{\ThinZZ}{\thth{\ZZ}}
\newcommand{\ThinG}{\thth{\Gsup}}
\newcommand{\ThinHF}{\thth{\HFsup}}
\newcommand{\ThinKh}{\thth{\Khsup}}
\makemirror{\mirrorThin}  {\Thin}
\makemirror{\mirrorThinZ} {\ThinZ} 
\makemirror{\mirrorThinZZ}{\ThinZZ} 
\makemirror{\mirrorThinG} {\ThinG} 
\makemirror{\mirrorThinHF}{\ThinHF} 
\makemirror{\mirrorThinKh}{\ThinKh} 

\makeBoundary{\BdryThin}  {\Thin}
\makeBoundary{\BdryThinZ} {\ThinZ} 
\makeBoundary{\BdryThinZZ}{\ThinZZ} 
\makeBoundary{\BdryThinG} {\ThinG} 
\makeBoundary{\BdryThinHF}{\ThinHF} 
\makeBoundary{\BdryThinKh}{\ThinKh} 
\makemirror{\mirrorBdryThin}  {\BdryThin}
\makemirror{\mirrorBdryThinZ} {\BdryThinZ} 
\makemirror{\mirrorBdryThinZZ}{\BdryThinZZ} 
\makemirror{\mirrorBdryThinG} {\BdryThinG} 
\makemirror{\mirrorBdryThinHF}{\BdryThinHF} 
\makemirror{\mirrorBdryThinKh}{\BdryThinKh} 

\makeInterior{\IntThin}  {\Thin}
\makeInterior{\IntThinZ} {\ThinZ} 
\makeInterior{\IntThinZZ}{\ThinZZ} 
\makeInterior{\IntThinG} {\ThinG} 
\makeInterior{\IntThinHF}{\ThinHF} 
\makeInterior{\IntThinKh}{\ThinKh} 
\makemirror{\mirrorIntThin}  {\IntThin} 
\makemirror{\mirrorIntThinZ} {\IntThinZ} 
\makemirror{\mirrorIntThinZZ}{\IntThinZZ} 
\makemirror{\mirrorIntThinG} {\IntThinG} 
\makemirror{\mirrorIntThinHF}{\IntThinHF} 
\makemirror{\mirrorIntThinKh}{\IntThinKh} 

\makeALink{\ALink}  {\Thin}
\makeALink{\ALinkHF}{\ThinHF} 
\makeALink{\ALinkKh}{\ThinKh} 
\makemirror{\mirrorALink}  {\ALink}
\makemirror{\mirrorALinkHF}{\ALinkHF} 
\makemirror{\mirrorALinkKh}{\ALinkKh} 

\makeBoundary{\BdryALink}  {\ALink}
\makeBoundary{\BdryALinkHF}{\ALinkHF} 
\makeBoundary{\BdryALinkKh}{\ALinkKh} 
\makemirror{\mirrorBdryALink}  {\BdryALink}
\makemirror{\mirrorBdryALinkHF}{\BdryALinkHF} 
\makemirror{\mirrorBdryALinkKh}{\BdryALinkKh} 

\makeInterior{\IntALink}  {\ALink}
\makeInterior{\IntALinkHF}{\ALinkHF} 
\makeInterior{\IntALinkKh}{\ALinkKh} 
\makemirror{\mirrorIntALink}  {\IntALink} 
\makemirror{\mirrorIntALinkHF}{\IntALinkHF} 
\makemirror{\mirrorIntALinkKh}{\IntALinkKh} 

\newcommand{\mirror}{\operatorname{m}} 

\newcommand{\FourPuncturedSphereKh}{S^2_{4,\ast}}

\newcommand{\QPI}{\operatorname{\mathbb{Q}P}^1}

\renewcommand{\c}{c}
\renewcommand{\a}{\mathbf{a}}
\newcommand{\tc}{\tilde{c}}
\newcommand{\s}{\mathbf{s}}
\newcommand{\ts}{\tilde{\s}}
\renewcommand{\r}{\mathbf{r}}
\newcommand{\tr}{\tilde{\r}}

%% file: PSTricks/input.tex
\newcommand{\vc}[1]{\vcenter{\hbox{#1}}}%
\newcommand{\mypic}[2]{%
	\newcommand{#2}{%
		\vc{%
			\includegraphics[page=#1]%
			{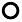}%
		}%
	}%
}%

\mypic{1}{\DotH}
\mypic{2}{\DotV}
\mypic{3}{\DotHblue}
\mypic{4}{\DotVblue}
\mypic{5}{\DotHred}
\mypic{6}{\DotVred}
\mypic{7}{\Li}
\mypic{8}{\Lo}
\mypic{9}{\LiRed}
\mypic{10}{\LoRed}
\mypic{11}{\LiBlue}
\mypic{12}{\LoBlue}
\mypic{13}{\Ni}
\mypic{14}{\No}
\mypic{15}{\KappaPairing}
\mypic{16}{\tanglepairing}
\mypic{17}{\ratixo}
\mypic{18}{\ratoxi}
\mypic{19}{\ratixi}
\mypic{20}{\ratixii}
\mypic{21}{\ratixiii}

%% file: computations/summary
\myrow{}{3_{1}}{-4}{1}{}\\
\myrow{}{4_{1}^{\star}}{0}{2}{}\\
\myrow{}{5_{1}}{-8}{2}{}\\
\myrow{}{5_{2}^{\star}}{-4}{3}{}\\
\myrow{}{6_{1}^{\star}}{0}{4}{}\\
\myrow{}{6_{2}^{\star}}{-4}{5}{}\\
\myrow{}{6_{3}^{\star}}{0}{6}{}\\
\myrow{}{7_{1}}{-12}{3}{}\\
\myrow{}{7_{2}^{\star}}{-4}{5}{}\\
\myrow{}{7_{3}^{\star}}{8}{6}{}\\
\myrow{}{7_{4}^{\star}}{4}{7}{}\\
\myrow{}{7_{5}^{\star}}{-8}{8}{}\\
\myrow{}{7_{6}^{\star}}{-4}{9}{}\\
\myrow{}{7_{7}^{\star}}{0}{10}{}\\
\myrow{}{8_{1}^{\star}}{0}{6}{}\\
\myrow{}{8_{2}^{\star}}{-8}{8}{}\\
\myrow{}{8_{3}^{\star}}{0}{8}{}\\
\myrow{}{8_{4}^{\star}}{-4}{9}{}\\
\myrow{}{8_{5}^{\star}}{8}{10}{}\\
\myrow{}{8_{6}^{\star}}{-4}{11}{}\\
\myrow{}{8_{7}^{\star}}{4}{11}{}\\
\myrow{}{8_{8}^{\star}}{0}{12}{}\\
\myrow{}{8_{9}^{\star}}{0}{12}{}\\
\myrow{}{8_{10}}{4}{13}{}\\
\myrow{}{8_{11}^{\star}}{-4}{13}{}\\
\myrow{}{8_{12}^{\star}}{0}{14}{}\\
\myrow{}{8_{13}^{\star}}{0}{14}{}\\
\myrow{}{8_{14}^{\star}}{4}{15}{}\\
\myrow{}{8_{15}^{\star}}{-8}{16}{}\\
\myrow{}{8_{16}}{-4}{17}{}\\
\myrow{}{8_{18}^{\star}}{0}{22}{}\\
\myrow{}{8_{19}}{8}{2}{}\\
\myrow{}{8_{20}}{0}{4}{}\\
\myrow{}{8_{21}^{\star}}{-4}{7}{}\\
\myrow{}{9_{1}}{-16}{4}{}\\
\myrow{}{9_{2}^{\star}}{-4}{7}{}\\
\myrow{}{9_{3}^{\star}}{12}{9}{}\\
\myrow{}{9_{4}^{\star}}{-8}{10}{}\\
\myrow{}{9_{5}^{\star}}{4}{11}{}\\
\myrow{}{9_{6}^{\star}}{-12}{13}{}\\
\myrow{}{9_{7}^{\star}}{-8}{14}{}\\
\myrow{}{9_{8}^{\star}}{-4}{15}{}\\
\myrow{}{9_{9}^{\star}}{-12}{15}{}\\
\myrow{}{9_{10}^{\star}}{8}{16}{}\\
\myrow{}{9_{11}^{\star}}{8}{16}{}\\
\myrow{}{9_{12}^{\star}}{-4}{17}{}\\
\myrow{}{9_{13}^{\star}}{8}{18}{}\\
\myrow{}{9_{14}^{\star}}{0}{18}{}\\
\myrow{}{9_{15}^{\star}}{4}{19}{}\\
\myrow{}{9_{16}^{\star}}{12}{19}{}\\
\myrow{}{9_{17}^{\star}}{-4}{19}{}\\
\myrow{}{9_{18}^{\star}}{-8}{20}{}\\
\myrow{}{9_{19}^{\star}}{0}{20}{}\\
\myrow{}{9_{20}^{\star}}{-8}{20}{}\\
\myrow{}{9_{21}^{\star}}{4}{21}{}\\
\myrow{}{9_{22}}{4}{21}{}\\
\myrow{}{9_{23}^{\star}}{-8}{22}{}\\
\myrow{}{9_{24}}{0}{22}{}\\
\myrow{}{9_{25}}{-4}{23}{}\\
\myrow{}{9_{26}^{\star}}{4}{23}{}\\
\myrow{}{9_{27}^{\star}}{0}{24}{}\\
\myrow{}{9_{28}^{\star}}{-4}{25}{}\\
\myrow{}{9_{29}}{-4}{25}{}\\
\myrow{}{9_{30}}{0}{26}{}\\
\myrow{}{9_{31}^{\star}}{4}{27}{}\\
\myrow{}{9_{34}}{0}{34}{}\\
\myrow{}{9_{35}^{\star}}{-4}{13}{}\\
\myrow{}{9_{36}}{8}{18}{}\\
\myrow{}{9_{37}^{\star}}{0}{22}{}\\
\myrow{}{9_{38}}{-8}{28}{}\\
\myrow{}{9_{39}}{4}{27}{}\\
\myrow{}{9_{40}^{\star}}{-4}{37}{}\\
\myrow{}{9_{41}}{0}{24}{}\\
\myrow{}{9_{42}}{0}{4}{}\\
\myrow{}{9_{43}}{8}{6}{}\\
\myrow{}{9_{44}}{0}{8}{}\\
\myrow{}{9_{45}}{-4}{11}{}\\
\myrow{\rowcolor{lightgray}}{9_{46}^{1}}{0}{4}{}\\
\myrow{\rowcolor{lightgray}}{9_{46}^{2}}{-4}{7}{}\\
\myrow{}{9_{47}}{4}{13}{}\\
\myrow{}{9_{48}^{\star}}{-4}{13}{}\\
\myrow{}{9_{49}}{8}{12}{}\\

%% file: computations/gradings
$
\begin{tikzpicture}[yscale=.25,xscale=.25]
\mynode{2}{4}{1}
\mygrid{2}{ 4}{ 2}{ 4}
\mylabel{2}{ 4}{ 2}{ 4}{ 3_{1}}
\end{tikzpicture}
$
\\
$
\begin{tikzpicture}[yscale=.25,xscale=.25]
\mynode{6}{4}{1}
 \mynode{2}{5}{1}
\mygrid{2}{ 4}{ 6}{ 5}
\mylabel{2}{ 4}{ 6}{ 5}{ 4_{1}^{1}}
\end{tikzpicture}
$
\\
$
\begin{tikzpicture}[yscale=.25,xscale=.25]
\mynode{-3}{-6}{1}
 \mynode{-7}{-5}{1}
\mygrid{-7}{ -6}{ -3}{ -5}
\mylabel{-7}{ -6}{ -3}{ -5}{ 4_{1}^{2}}
\end{tikzpicture}
$
\\
$
\begin{tikzpicture}[yscale=.25,xscale=.25]
\mynode{-2}{4}{1}
 \mynode{-1}{4}{1}
\mygrid{-2}{ 4}{ -1}{ 4}
\mylabel{-2}{ 4}{ -1}{ 4}{ 5_{1}}
\end{tikzpicture}
$
\\
$
\begin{tikzpicture}[yscale=.25,xscale=.25]
\mynode{10}{8}{1}
 \mynode{6}{9}{1}
 \mynode{2}{10}{1}
\mygrid{2}{ 8}{ 10}{ 10}
\mylabel{2}{ 8}{ 10}{ 10}{ 5_{2}^{1}}
\end{tikzpicture}
$
\\
$
\begin{tikzpicture}[yscale=.25,xscale=.25]
\mynode{-3}{-1}{1}
 \mynode{-7}{0}{1}
 \mynode{-6}{0}{1}
\mygrid{-7}{ -1}{ -3}{ 0}
\mylabel{-7}{ -1}{ -3}{ 0}{ 5_{2}^{2}}
\end{tikzpicture}
$
\\
$
\begin{tikzpicture}[yscale=.25,xscale=.25]
\mynode{14}{8}{1}
 \mynode{10}{9}{1}
 \mynode{6}{10}{1}
 \mynode{2}{11}{1}
\mygrid{2}{ 8}{ 14}{ 11}
\mylabel{2}{ 8}{ 14}{ 11}{ 6_{1}^{1}}
\end{tikzpicture}
$
\\
$
\begin{tikzpicture}[yscale=.25,xscale=.25]
\mynode{4}{-2}{1}
 \mynode{5}{-2}{1}
 \mynode{0}{-1}{1}
 \mynode{1}{-1}{1}
\mygrid{0}{ -2}{ 5}{ -1}
\mylabel{0}{ -2}{ 5}{ -1}{ 6_{1}^{2}}
\end{tikzpicture}
$
\\
$
\begin{tikzpicture}[yscale=.25,xscale=.25]
\mynode{2}{4}{2}
 \mynode{3}{4}{1}
 \mynode{-2}{5}{1}
 \mynode{-1}{5}{1}
\mygrid{-2}{ 4}{ 3}{ 5}
\mylabel{-2}{ 4}{ 3}{ 5}{ 6_{2}^{1}}
\end{tikzpicture}
$
\\
$
\begin{tikzpicture}[yscale=.25,xscale=.25]
\mynode{9}{2}{1}
 \mynode{10}{2}{1}
 \mynode{5}{3}{1}
 \mynode{6}{3}{1}
 \mynode{2}{4}{1}
\mygrid{2}{ 2}{ 10}{ 4}
\mylabel{2}{ 2}{ 10}{ 4}{ 6_{2}^{2}}
\end{tikzpicture}
$
\\
$
\begin{tikzpicture}[yscale=.25,xscale=.25]
\mynode{-3}{-6}{1}
 \mynode{-7}{-5}{2}
 \mynode{-6}{-5}{1}
 \mynode{-11}{-4}{1}
 \mynode{-10}{-4}{1}
\mygrid{-11}{ -6}{ -3}{ -4}
\mylabel{-11}{ -6}{ -3}{ -4}{ 6_{3}^{1}}
\end{tikzpicture}
$
\\
$
\begin{tikzpicture}[yscale=.25,xscale=.25]
\mynode{9}{3}{1}
 \mynode{10}{3}{1}
 \mynode{5}{4}{1}
 \mynode{6}{4}{2}
 \mynode{2}{5}{1}
\mygrid{2}{ 3}{ 10}{ 5}
\mylabel{2}{ 3}{ 10}{ 5}{ 6_{3}^{2}}
\end{tikzpicture}
$
\\
$
\begin{tikzpicture}[yscale=.25,xscale=.25]
\mynode{2}{8}{1}
 \mynode{3}{8}{1}
 \mynode{4}{8}{1}
\mygrid{2}{ 8}{ 4}{ 8}
\mylabel{2}{ 8}{ 4}{ 8}{ 7_{1}}
\end{tikzpicture}
$
\\
$
\begin{tikzpicture}[yscale=.25,xscale=.25]
\mynode{18}{12}{1}
 \mynode{14}{13}{1}
 \mynode{10}{14}{1}
 \mynode{6}{15}{1}
 \mynode{2}{16}{1}
\mygrid{2}{ 12}{ 18}{ 16}
\mylabel{2}{ 12}{ 18}{ 16}{ 7_{2}^{1}}
\end{tikzpicture}
$
\\
$
\begin{tikzpicture}[yscale=.25,xscale=.25]
\mynode{4}{3}{1}
 \mynode{5}{3}{1}
 \mynode{0}{4}{1}
 \mynode{1}{4}{1}
 \mynode{2}{4}{1}
\mygrid{0}{ 3}{ 5}{ 4}
\mylabel{0}{ 3}{ 5}{ 4}{ 7_{2}^{2}}
\end{tikzpicture}
$
\\
$
\begin{tikzpicture}[yscale=.25,xscale=.25]
\mynode{8}{-7}{1}
 \mynode{9}{-7}{1}
 \mynode{4}{-6}{1}
 \mynode{5}{-6}{1}
 \mynode{0}{-5}{1}
 \mynode{1}{-5}{1}
\mygrid{0}{ -7}{ 9}{ -5}
\mylabel{0}{ -7}{ 9}{ -5}{ 7_{3}^{1}}
\end{tikzpicture}
$
\\
$
\begin{tikzpicture}[yscale=.25,xscale=.25]
\mynode{0}{-5}{1}
 \mynode{1}{-5}{2}
 \mynode{2}{-5}{1}
 \mynode{-3}{-4}{1}
 \mynode{-2}{-4}{1}
\mygrid{-3}{ -5}{ 2}{ -4}
\mylabel{-3}{ -5}{ 2}{ -4}{ 7_{3}^{2}}
\end{tikzpicture}
$
\\
$
\begin{tikzpicture}[yscale=.25,xscale=.25]
\mynode{13}{-3}{1}
 \mynode{14}{-3}{1}
 \mynode{9}{-2}{1}
 \mynode{10}{-2}{1}
 \mynode{5}{-1}{1}
 \mynode{6}{-1}{1}
 \mynode{2}{0}{1}
\mygrid{2}{ -3}{ 14}{ 0}
\mylabel{2}{ -3}{ 14}{ 0}{ 7_{4}^{1}}
\end{tikzpicture}
$
\\
$
\begin{tikzpicture}[yscale=.25,xscale=.25]
\mynode{-3}{-11}{1}
 \mynode{-7}{-10}{1}
 \mynode{-11}{-9}{1}
 \mynode{-10}{-9}{1}
 \mynode{-9}{-9}{1}
 \mynode{-14}{-8}{1}
 \mynode{-13}{-8}{1}
\mygrid{-14}{ -11}{ -3}{ -8}
\mylabel{-14}{ -11}{ -3}{ -8}{ 7_{4}^{2}}
\end{tikzpicture}
$
\\
$
\begin{tikzpicture}[yscale=.25,xscale=.25]
\mynode{6}{8}{2}
 \mynode{7}{8}{1}
 \mynode{2}{9}{2}
 \mynode{3}{9}{1}
 \mynode{-2}{10}{1}
 \mynode{-1}{10}{1}
\mygrid{-2}{ 8}{ 7}{ 10}
\mylabel{-2}{ 8}{ 7}{ 10}{ 7_{5}^{1}}
\end{tikzpicture}
$
\\
$
\begin{tikzpicture}[yscale=.25,xscale=.25]
\mynode{9}{7}{1}
 \mynode{10}{7}{1}
 \mynode{5}{8}{1}
 \mynode{6}{8}{3}
 \mynode{7}{8}{1}
 \mynode{2}{9}{1}
\mygrid{2}{ 7}{ 10}{ 9}
\mylabel{2}{ 7}{ 10}{ 9}{ 7_{5}^{2}}
\end{tikzpicture}
$
\\
$
\begin{tikzpicture}[yscale=.25,xscale=.25]
\mynode{2}{4}{2}
 \mynode{3}{4}{1}
 \mynode{-2}{5}{2}
 \mynode{-1}{5}{2}
 \mynode{-6}{6}{1}
 \mynode{-5}{6}{1}
\mygrid{-6}{ 4}{ 3}{ 6}
\mylabel{-6}{ 4}{ 3}{ 6}{ 7_{6}^{1}}
\end{tikzpicture}
$
\\
$
\begin{tikzpicture}[yscale=.25,xscale=.25]
\mynode{-3}{-7}{1}
 \mynode{-7}{-6}{2}
 \mynode{-6}{-6}{1}
 \mynode{-11}{-5}{2}
 \mynode{-10}{-5}{1}
 \mynode{-15}{-4}{1}
 \mynode{-14}{-4}{1}
\mygrid{-15}{ -7}{ -3}{ -4}
\mylabel{-15}{ -7}{ -3}{ -4}{ 7_{6}^{2}}
\end{tikzpicture}
$
\\
$
\begin{tikzpicture}[yscale=.25,xscale=.25]
\mynode{-3}{-6}{1}
 \mynode{-7}{-5}{2}
 \mynode{-6}{-5}{1}
 \mynode{-11}{-4}{2}
 \mynode{-10}{-4}{2}
 \mynode{-15}{-3}{1}
 \mynode{-14}{-3}{1}
\mygrid{-15}{ -6}{ -3}{ -3}
\mylabel{-15}{ -6}{ -3}{ -3}{ 7_{7}^{1}}
\end{tikzpicture}
$
\\
$
\begin{tikzpicture}[yscale=.25,xscale=.25]
\mynode{12}{2}{1}
 \mynode{13}{2}{1}
 \mynode{8}{3}{1}
 \mynode{9}{3}{2}
 \mynode{10}{3}{1}
 \mynode{5}{4}{1}
 \mynode{6}{4}{2}
 \mynode{2}{5}{1}
\mygrid{2}{ 2}{ 13}{ 5}
\mylabel{2}{ 2}{ 13}{ 5}{ 7_{7}^{2}}
\end{tikzpicture}
$
\\
$
\begin{tikzpicture}[yscale=.25,xscale=.25]
\mynode{22}{12}{1}
 \mynode{18}{13}{1}
 \mynode{14}{14}{1}
 \mynode{10}{15}{1}
 \mynode{6}{16}{1}
 \mynode{2}{17}{1}
\mygrid{2}{ 12}{ 22}{ 17}
\mylabel{2}{ 12}{ 22}{ 17}{ 8_{1}^{1}}
\end{tikzpicture}
$
\\
$
\begin{tikzpicture}[yscale=.25,xscale=.25]
\mynode{-5}{-6}{1}
 \mynode{-4}{-6}{1}
 \mynode{-3}{-6}{1}
 \mynode{-9}{-5}{1}
 \mynode{-8}{-5}{1}
 \mynode{-7}{-5}{1}
\mygrid{-9}{ -6}{ -3}{ -5}
\mylabel{-9}{ -6}{ -3}{ -5}{ 8_{1}^{2}}
\end{tikzpicture}
$
\\
$
\begin{tikzpicture}[yscale=.25,xscale=.25]
\mynode{6}{8}{2}
 \mynode{7}{8}{2}
 \mynode{8}{8}{1}
 \mynode{2}{9}{1}
 \mynode{3}{9}{1}
 \mynode{4}{9}{1}
\mygrid{2}{ 8}{ 8}{ 9}
\mylabel{2}{ 8}{ 8}{ 9}{ 8_{2}^{1}}
\end{tikzpicture}
$
\\
$
\begin{tikzpicture}[yscale=.25,xscale=.25]
\mynode{-3}{-2}{1}
 \mynode{-2}{-2}{1}
 \mynode{-1}{-2}{1}
 \mynode{-7}{-1}{1}
 \mynode{-6}{-1}{1}
 \mynode{-5}{-1}{1}
 \mynode{-10}{0}{1}
 \mynode{-9}{0}{1}
\mygrid{-10}{ -2}{ -1}{ 0}
\mylabel{-10}{ -2}{ -1}{ 0}{ 8_{2}^{2}}
\end{tikzpicture}
$
\\
$
\begin{tikzpicture}[yscale=.25,xscale=.25]
\mynode{-2}{0}{1}
 \mynode{-1}{0}{1}
 \mynode{-6}{1}{1}
 \mynode{-5}{1}{1}
 \mynode{-10}{2}{1}
 \mynode{-9}{2}{1}
 \mynode{-14}{3}{1}
 \mynode{-13}{3}{1}
\mygrid{-14}{ 0}{ -1}{ 3}
\mylabel{-14}{ 0}{ -1}{ 3}{ 8_{3}^{1}}
\end{tikzpicture}
$
\\
$
\begin{tikzpicture}[yscale=.25,xscale=.25]
\mynode{12}{-4}{1}
 \mynode{13}{-4}{1}
 \mynode{8}{-3}{1}
 \mynode{9}{-3}{1}
 \mynode{4}{-2}{1}
 \mynode{5}{-2}{1}
 \mynode{0}{-1}{1}
 \mynode{1}{-1}{1}
\mygrid{0}{ -4}{ 13}{ -1}
\mylabel{0}{ -4}{ 13}{ -1}{ 8_{3}^{2}}
\end{tikzpicture}
$
\\
$
\begin{tikzpicture}[yscale=.25,xscale=.25]
\mynode{17}{0}{1}
 \mynode{18}{0}{1}
 \mynode{13}{1}{1}
 \mynode{14}{1}{1}
 \mynode{9}{2}{1}
 \mynode{10}{2}{1}
 \mynode{5}{3}{1}
 \mynode{6}{3}{1}
 \mynode{2}{4}{1}
\mygrid{2}{ 0}{ 18}{ 4}
\mylabel{2}{ 0}{ 18}{ 4}{ 8_{4}^{1}}
\end{tikzpicture}
$
\\
$
\begin{tikzpicture}[yscale=.25,xscale=.25]
\mynode{2}{4}{2}
 \mynode{3}{4}{2}
 \mynode{4}{4}{1}
 \mynode{-2}{5}{1}
 \mynode{-1}{5}{2}
 \mynode{0}{5}{1}
\mygrid{-2}{ 4}{ 4}{ 5}
\mylabel{-2}{ 4}{ 4}{ 5}{ 8_{4}^{2}}
\end{tikzpicture}
$
\\
$
\begin{tikzpicture}[yscale=.25,xscale=.25]
\mynode{3}{-6}{1}
 \mynode{4}{-6}{2}
 \mynode{5}{-6}{1}
 \mynode{-1}{-5}{1}
 \mynode{0}{-5}{3}
 \mynode{1}{-5}{2}
\mygrid{-1}{ -6}{ 5}{ -5}
\mylabel{-1}{ -6}{ 5}{ -5}{ 8_{5}^{1}}
\end{tikzpicture}
$
\\
$
\begin{tikzpicture}[yscale=.25,xscale=.25]
\mynode{16}{3}{1}
 \mynode{17}{3}{1}
 \mynode{12}{4}{1}
 \mynode{13}{4}{1}
 \mynode{14}{4}{1}
 \mynode{8}{5}{1}
 \mynode{9}{5}{1}
 \mynode{10}{5}{1}
 \mynode{6}{6}{1}
 \mynode{2}{7}{1}
\mygrid{2}{ 3}{ 17}{ 7}
\mylabel{2}{ 3}{ 17}{ 7}{ 8_{5}^{2}}
\end{tikzpicture}
$
\\
$
\begin{tikzpicture}[yscale=.25,xscale=.25]
\mynode{0}{-2}{1}
 \mynode{1}{-2}{2}
 \mynode{2}{-2}{1}
 \mynode{-4}{-1}{1}
 \mynode{-3}{-1}{3}
 \mynode{-2}{-1}{1}
 \mynode{-7}{0}{1}
 \mynode{-6}{0}{1}
\mygrid{-7}{ -2}{ 2}{ 0}
\mylabel{-7}{ -2}{ 2}{ 0}{ 8_{6}^{1}}
\end{tikzpicture}
$
\\
$
\begin{tikzpicture}[yscale=.25,xscale=.25]
\mynode{10}{8}{2}
 \mynode{11}{8}{1}
 \mynode{6}{9}{2}
 \mynode{7}{9}{1}
 \mynode{2}{10}{2}
 \mynode{3}{10}{1}
 \mynode{-2}{11}{1}
 \mynode{-1}{11}{1}
\mygrid{-2}{ 8}{ 11}{ 11}
\mylabel{-2}{ 8}{ 11}{ 11}{ 8_{6}^{2}}
\end{tikzpicture}
$
\\
$
\begin{tikzpicture}[yscale=.25,xscale=.25]
\mynode{8}{-2}{1}
 \mynode{9}{-2}{1}
 \mynode{4}{-1}{2}
 \mynode{5}{-1}{3}
 \mynode{6}{-1}{1}
 \mynode{0}{0}{1}
 \mynode{1}{0}{1}
 \mynode{2}{0}{1}
\mygrid{0}{ -2}{ 9}{ 0}
\mylabel{0}{ -2}{ 9}{ 0}{ 8_{7}^{1}}
\end{tikzpicture}
$
\\
$
\begin{tikzpicture}[yscale=.25,xscale=.25]
\mynode{4}{-1}{1}
 \mynode{5}{-1}{2}
 \mynode{6}{-1}{1}
 \mynode{0}{0}{1}
 \mynode{1}{0}{2}
 \mynode{2}{0}{2}
 \mynode{-3}{1}{1}
 \mynode{-2}{1}{1}
\mygrid{-3}{ -1}{ 6}{ 1}
\mylabel{-3}{ -1}{ 6}{ 1}{ 8_{7}^{2}}
\end{tikzpicture}
$
\\
$
\begin{tikzpicture}[yscale=.25,xscale=.25]
\mynode{1}{-1}{1}
 \mynode{2}{-1}{2}
 \mynode{3}{-1}{1}
 \mynode{-3}{0}{1}
 \mynode{-2}{0}{3}
 \mynode{-1}{0}{2}
 \mynode{-6}{1}{1}
 \mynode{-5}{1}{1}
\mygrid{-6}{ -1}{ 3}{ 1}
\mylabel{-6}{ -1}{ 3}{ 1}{ 8_{8}^{1}}
\end{tikzpicture}
$
\\
$
\begin{tikzpicture}[yscale=.25,xscale=.25]
\mynode{-3}{-12}{1}
 \mynode{-7}{-11}{2}
 \mynode{-6}{-11}{1}
 \mynode{-11}{-10}{2}
 \mynode{-10}{-10}{1}
 \mynode{-15}{-9}{2}
 \mynode{-14}{-9}{1}
 \mynode{-19}{-8}{1}
 \mynode{-18}{-8}{1}
\mygrid{-19}{ -12}{ -3}{ -8}
\mylabel{-19}{ -12}{ -3}{ -8}{ 8_{8}^{2}}
\end{tikzpicture}
$
\\
$
\begin{tikzpicture}[yscale=.25,xscale=.25]
\mynode{1}{-1}{1}
 \mynode{2}{-1}{1}
 \mynode{-3}{0}{2}
 \mynode{-2}{0}{3}
 \mynode{-1}{0}{1}
 \mynode{-7}{1}{1}
 \mynode{-6}{1}{2}
 \mynode{-5}{1}{1}
\mygrid{-7}{ -1}{ 2}{ 1}
\mylabel{-7}{ -1}{ 2}{ 1}{ 8_{9}^{1}}
\end{tikzpicture}
$
\\
$
\begin{tikzpicture}[yscale=.25,xscale=.25]
\mynode{4}{-2}{1}
 \mynode{5}{-2}{2}
 \mynode{6}{-2}{1}
 \mynode{0}{-1}{1}
 \mynode{1}{-1}{3}
 \mynode{2}{-1}{2}
 \mynode{-3}{0}{1}
 \mynode{-2}{0}{1}
\mygrid{-3}{ -2}{ 6}{ 0}
\mylabel{-3}{ -2}{ 6}{ 0}{ 8_{9}^{2}}
\end{tikzpicture}
$
\\
$
\begin{tikzpicture}[yscale=.25,xscale=.25]
\mynode{0}{-6}{1}
 \mynode{1}{-6}{1}
 \mynode{-4}{-5}{2}
 \mynode{-3}{-5}{4}
 \mynode{-2}{-5}{1}
 \mynode{-8}{-4}{1}
 \mynode{-7}{-4}{2}
 \mynode{-6}{-4}{1}
\mygrid{-8}{ -6}{ 1}{ -4}
\mylabel{-8}{ -6}{ 1}{ -4}{ 8_{10}}
\end{tikzpicture}
$
\\
$
\begin{tikzpicture}[yscale=.25,xscale=.25]
\mynode{2}{4}{2}
 \mynode{3}{4}{1}
 \mynode{-2}{5}{2}
 \mynode{-1}{5}{2}
 \mynode{-6}{6}{2}
 \mynode{-5}{6}{2}
 \mynode{-10}{7}{1}
 \mynode{-9}{7}{1}
\mygrid{-10}{ 4}{ 3}{ 7}
\mylabel{-10}{ 4}{ 3}{ 7}{ 8_{11}^{1}}
\end{tikzpicture}
$
\\
$
\begin{tikzpicture}[yscale=.25,xscale=.25]
\mynode{-3}{-7}{1}
 \mynode{-8}{-6}{1}
 \mynode{-7}{-6}{3}
 \mynode{-6}{-6}{1}
 \mynode{-12}{-5}{1}
 \mynode{-11}{-5}{3}
 \mynode{-10}{-5}{1}
 \mynode{-15}{-4}{1}
 \mynode{-14}{-4}{1}
\mygrid{-15}{ -7}{ -3}{ -4}
\mylabel{-15}{ -7}{ -3}{ -4}{ 8_{11}^{2}}
\end{tikzpicture}
$
\\
$
\begin{tikzpicture}[yscale=.25,xscale=.25]
\mynode{6}{4}{2}
 \mynode{7}{4}{1}
 \mynode{2}{5}{3}
 \mynode{3}{5}{2}
 \mynode{-2}{6}{2}
 \mynode{-1}{6}{2}
 \mynode{-6}{7}{1}
 \mynode{-5}{7}{1}
\mygrid{-6}{ 4}{ 7}{ 7}
\mylabel{-6}{ 4}{ 7}{ 7}{ 8_{12}^{1}}
\end{tikzpicture}
$
\\
$
\begin{tikzpicture}[yscale=.25,xscale=.25]
\mynode{4}{-8}{1}
 \mynode{5}{-8}{1}
 \mynode{0}{-7}{2}
 \mynode{1}{-7}{2}
 \mynode{-4}{-6}{2}
 \mynode{-3}{-6}{3}
 \mynode{-8}{-5}{1}
 \mynode{-7}{-5}{2}
\mygrid{-8}{ -8}{ 5}{ -5}
\mylabel{-8}{ -8}{ 5}{ -5}{ 8_{12}^{2}}
\end{tikzpicture}
$
\\
$
\begin{tikzpicture}[yscale=.25,xscale=.25]
\mynode{-3}{-6}{1}
 \mynode{-7}{-5}{2}
 \mynode{-6}{-5}{1}
 \mynode{-11}{-4}{2}
 \mynode{-10}{-4}{2}
 \mynode{-15}{-3}{2}
 \mynode{-14}{-3}{2}
 \mynode{-19}{-2}{1}
 \mynode{-18}{-2}{1}
\mygrid{-19}{ -6}{ -3}{ -2}
\mylabel{-19}{ -6}{ -3}{ -2}{ 8_{13}^{1}}
\end{tikzpicture}
$
\\
$
\begin{tikzpicture}[yscale=.25,xscale=.25]
\mynode{-3}{-6}{1}
 \mynode{-8}{-5}{1}
 \mynode{-7}{-5}{3}
 \mynode{-6}{-5}{1}
 \mynode{-12}{-4}{1}
 \mynode{-11}{-4}{3}
 \mynode{-10}{-4}{2}
 \mynode{-15}{-3}{1}
 \mynode{-14}{-3}{1}
\mygrid{-15}{ -6}{ -3}{ -3}
\mylabel{-15}{ -6}{ -3}{ -3}{ 8_{13}^{2}}
\end{tikzpicture}
$
\\
$
\begin{tikzpicture}[yscale=.25,xscale=.25]
\mynode{-3}{-11}{1}
 \mynode{-8}{-10}{1}
 \mynode{-7}{-10}{3}
 \mynode{-6}{-10}{1}
 \mynode{-12}{-9}{1}
 \mynode{-11}{-9}{4}
 \mynode{-10}{-9}{2}
 \mynode{-15}{-8}{1}
 \mynode{-14}{-8}{1}
\mygrid{-15}{ -11}{ -3}{ -8}
\mylabel{-15}{ -11}{ -3}{ -8}{ 8_{14}^{1}}
\end{tikzpicture}
$
\\
$
\begin{tikzpicture}[yscale=.25,xscale=.25]
\mynode{-3}{-11}{1}
 \mynode{-7}{-10}{2}
 \mynode{-6}{-10}{1}
 \mynode{-11}{-9}{3}
 \mynode{-10}{-9}{2}
 \mynode{-15}{-8}{2}
 \mynode{-14}{-8}{2}
 \mynode{-19}{-7}{1}
 \mynode{-18}{-7}{1}
\mygrid{-19}{ -11}{ -3}{ -7}
\mylabel{-19}{ -11}{ -3}{ -7}{ 8_{14}^{2}}
\end{tikzpicture}
$
\\
$
\begin{tikzpicture}[yscale=.25,xscale=.25]
\mynode{-3}{-2}{1}
 \mynode{-7}{-1}{3}
 \mynode{-6}{-1}{2}
 \mynode{-11}{0}{3}
 \mynode{-10}{0}{4}
 \mynode{-9}{0}{1}
 \mynode{-15}{1}{1}
 \mynode{-14}{1}{1}
\mygrid{-15}{ -2}{ -3}{ 1}
\mylabel{-15}{ -2}{ -3}{ 1}{ 8_{15}^{1}}
\end{tikzpicture}
$
\\
$
\begin{tikzpicture}[yscale=.25,xscale=.25]
\mynode{-3}{-2}{1}
 \mynode{-7}{-1}{2}
 \mynode{-6}{-1}{1}
 \mynode{-11}{0}{2}
 \mynode{-10}{0}{3}
 \mynode{-9}{0}{1}
 \mynode{-15}{1}{1}
 \mynode{-14}{1}{2}
 \mynode{-13}{1}{1}
 \mynode{-18}{2}{1}
 \mynode{-17}{2}{1}
\mygrid{-18}{ -2}{ -3}{ 2}
\mylabel{-18}{ -2}{ -3}{ 2}{ 8_{15}^{2}}
\end{tikzpicture}
$
\\
$
\begin{tikzpicture}[yscale=.25,xscale=.25]
\mynode{-3}{-7}{1}
 \mynode{-7}{-6}{2}
 \mynode{-6}{-6}{1}
 \mynode{-11}{-5}{2}
 \mynode{-10}{-5}{2}
 \mynode{-9}{-5}{1}
 \mynode{-15}{-4}{1}
 \mynode{-14}{-4}{3}
 \mynode{-13}{-4}{2}
 \mynode{-18}{-3}{1}
 \mynode{-17}{-3}{1}
\mygrid{-18}{ -7}{ -3}{ -3}
\mylabel{-18}{ -7}{ -3}{ -3}{ 8_{16}}
\end{tikzpicture}
$
\\
$
\begin{tikzpicture}[yscale=.25,xscale=.25]
\mynode{-3}{-6}{1}
 \mynode{-7}{-5}{3}
 \mynode{-6}{-5}{2}
 \mynode{-11}{-4}{3}
 \mynode{-10}{-4}{4}
 \mynode{-9}{-4}{1}
 \mynode{-15}{-3}{1}
 \mynode{-14}{-3}{3}
 \mynode{-13}{-3}{2}
 \mynode{-18}{-2}{1}
 \mynode{-17}{-2}{1}
\mygrid{-18}{ -6}{ -3}{ -2}
\mylabel{-18}{ -6}{ -3}{ -2}{ 8_{18}^{1}}
\end{tikzpicture}
$
\\
$
\begin{tikzpicture}[yscale=.25,xscale=.25]
\mynode{16}{1}{1}
 \mynode{17}{1}{1}
 \mynode{12}{2}{2}
 \mynode{13}{2}{3}
 \mynode{14}{2}{1}
 \mynode{8}{3}{1}
 \mynode{9}{3}{4}
 \mynode{10}{3}{3}
 \mynode{5}{4}{2}
 \mynode{6}{4}{3}
 \mynode{2}{5}{1}
\mygrid{2}{ 1}{ 17}{ 5}
\mylabel{2}{ 1}{ 17}{ 5}{ 8_{18}^{2}}
\end{tikzpicture}
$
\\
$
\begin{tikzpicture}[yscale=.25,xscale=.25]
\mynode{3}{-5}{1}
 \mynode{4}{-5}{1}
\mygrid{3}{ -5}{ 4}{ -5}
\mylabel{3}{ -5}{ 4}{ -5}{ 8_{19}}
\end{tikzpicture}
$
\\
$
\begin{tikzpicture}[yscale=.25,xscale=.25]
\mynode{0}{-1}{1}
 \mynode{1}{-1}{1}
 \mynode{-4}{0}{1}
 \mynode{-3}{0}{1}
\mygrid{-4}{ -1}{ 1}{ 0}
\mylabel{-4}{ -1}{ 1}{ 0}{ 8_{20}}
\end{tikzpicture}
$
\\
$
\begin{tikzpicture}[yscale=.25,xscale=.25]
\mynode{-3}{-1}{1}
 \mynode{-7}{0}{2}
 \mynode{-6}{0}{2}
 \mynode{-11}{1}{1}
 \mynode{-10}{1}{1}
\mygrid{-11}{ -1}{ -3}{ 1}
\mylabel{-11}{ -1}{ -3}{ 1}{ 8_{21}^{1}}
\end{tikzpicture}
$
\\
$
\begin{tikzpicture}[yscale=.25,xscale=.25]
\mynode{0}{-2}{1}
 \mynode{1}{-2}{1}
 \mynode{-4}{-1}{1}
 \mynode{-3}{-1}{2}
 \mynode{-7}{0}{1}
 \mynode{-6}{0}{1}
\mygrid{-7}{ -2}{ 1}{ 0}
\mylabel{-7}{ -2}{ 1}{ 0}{ 8_{21}^{2}}
\end{tikzpicture}
$
\\
$
\begin{tikzpicture}[yscale=.25,xscale=.25]
\mynode{-2}{8}{1}
 \mynode{-1}{8}{1}
 \mynode{0}{8}{1}
 \mynode{1}{8}{1}
\mygrid{-2}{ 8}{ 1}{ 8}
\mylabel{-2}{ 8}{ 1}{ 8}{ 9_{1}}
\end{tikzpicture}
$
\\
$
\begin{tikzpicture}[yscale=.25,xscale=.25]
\mynode{26}{16}{1}
 \mynode{22}{17}{1}
 \mynode{18}{18}{1}
 \mynode{14}{19}{1}
 \mynode{10}{20}{1}
 \mynode{6}{21}{1}
 \mynode{2}{22}{1}
\mygrid{2}{ 16}{ 26}{ 22}
\mylabel{2}{ 16}{ 26}{ 22}{ 9_{2}^{1}}
\end{tikzpicture}
$
\\
$
\begin{tikzpicture}[yscale=.25,xscale=.25]
\mynode{-5}{-1}{1}
 \mynode{-4}{-1}{1}
 \mynode{-3}{-1}{1}
 \mynode{-9}{0}{1}
 \mynode{-8}{0}{1}
 \mynode{-7}{0}{1}
 \mynode{-6}{0}{1}
\mygrid{-9}{ -1}{ -3}{ 0}
\mylabel{-9}{ -1}{ -3}{ 0}{ 9_{2}^{2}}
\end{tikzpicture}
$
\\
$
\begin{tikzpicture}[yscale=.25,xscale=.25]
\mynode{-5}{-15}{1}
 \mynode{-4}{-15}{1}
 \mynode{-3}{-15}{1}
 \mynode{-9}{-14}{1}
 \mynode{-8}{-14}{1}
 \mynode{-7}{-14}{1}
 \mynode{-13}{-13}{1}
 \mynode{-12}{-13}{1}
 \mynode{-11}{-13}{1}
\mygrid{-13}{ -15}{ -3}{ -13}
\mylabel{-13}{ -15}{ -3}{ -13}{ 9_{3}^{1}}
\end{tikzpicture}
$
\\
$
\begin{tikzpicture}[yscale=.25,xscale=.25]
\mynode{3}{-5}{1}
 \mynode{4}{-5}{2}
 \mynode{5}{-5}{2}
 \mynode{6}{-5}{1}
 \mynode{0}{-4}{1}
 \mynode{1}{-4}{1}
 \mynode{2}{-4}{1}
\mygrid{0}{ -5}{ 6}{ -4}
\mylabel{0}{ -5}{ 6}{ -4}{ 9_{3}^{2}}
\end{tikzpicture}
$
\\
$
\begin{tikzpicture}[yscale=.25,xscale=.25]
\mynode{-2}{4}{1}
 \mynode{-1}{4}{1}
 \mynode{-6}{5}{1}
 \mynode{-5}{5}{1}
 \mynode{-10}{6}{1}
 \mynode{-9}{6}{1}
 \mynode{-14}{7}{1}
 \mynode{-13}{7}{1}
 \mynode{-18}{8}{1}
 \mynode{-17}{8}{1}
\mygrid{-18}{ 4}{ -1}{ 8}
\mylabel{-18}{ 4}{ -1}{ 8}{ 9_{4}^{1}}
\end{tikzpicture}
$
\\
$
\begin{tikzpicture}[yscale=.25,xscale=.25]
\mynode{0}{3}{1}
 \mynode{1}{3}{2}
 \mynode{2}{3}{1}
 \mynode{-4}{4}{1}
 \mynode{-3}{4}{2}
 \mynode{-2}{4}{2}
 \mynode{-1}{4}{1}
\mygrid{-4}{ 3}{ 2}{ 4}
\mylabel{-4}{ 3}{ 2}{ 4}{ 9_{4}^{2}}
\end{tikzpicture}
$
\\
$
\begin{tikzpicture}[yscale=.25,xscale=.25]
\mynode{21}{-5}{1}
 \mynode{22}{-5}{1}
 \mynode{17}{-4}{1}
 \mynode{18}{-4}{1}
 \mynode{13}{-3}{1}
 \mynode{14}{-3}{1}
 \mynode{9}{-2}{1}
 \mynode{10}{-2}{1}
 \mynode{5}{-1}{1}
 \mynode{6}{-1}{1}
 \mynode{2}{0}{1}
\mygrid{2}{ -5}{ 22}{ 0}
\mylabel{2}{ -5}{ 22}{ 0}{ 9_{5}^{1}}
\end{tikzpicture}
$
\\
$
\begin{tikzpicture}[yscale=.25,xscale=.25]
\mynode{-3}{-11}{1}
 \mynode{-2}{-11}{1}
 \mynode{-1}{-11}{1}
 \mynode{-7}{-10}{1}
 \mynode{-6}{-10}{1}
 \mynode{-5}{-10}{1}
 \mynode{-11}{-9}{1}
 \mynode{-10}{-9}{1}
 \mynode{-9}{-9}{1}
 \mynode{-14}{-8}{1}
 \mynode{-13}{-8}{1}
\mygrid{-14}{ -11}{ -1}{ -8}
\mylabel{-14}{ -11}{ -1}{ -8}{ 9_{5}^{2}}
\end{tikzpicture}
$
\\
$
\begin{tikzpicture}[yscale=.25,xscale=.25]
\mynode{10}{12}{2}
 \mynode{11}{12}{2}
 \mynode{12}{12}{1}
 \mynode{6}{13}{2}
 \mynode{7}{13}{2}
 \mynode{8}{13}{1}
 \mynode{2}{14}{1}
 \mynode{3}{14}{1}
 \mynode{4}{14}{1}
\mygrid{2}{ 12}{ 12}{ 14}
\mylabel{2}{ 12}{ 12}{ 14}{ 9_{6}^{1}}
\end{tikzpicture}
$
\\
$
\begin{tikzpicture}[yscale=.25,xscale=.25]
\mynode{-3}{3}{1}
 \mynode{-2}{3}{1}
 \mynode{-1}{3}{1}
 \mynode{-7}{4}{1}
 \mynode{-6}{4}{3}
 \mynode{-5}{4}{3}
 \mynode{-4}{4}{1}
 \mynode{-10}{5}{1}
 \mynode{-9}{5}{1}
\mygrid{-10}{ 3}{ -1}{ 5}
\mylabel{-10}{ 3}{ -1}{ 5}{ 9_{6}^{2}}
\end{tikzpicture}
$
\\
$
\begin{tikzpicture}[yscale=.25,xscale=.25]
\mynode{14}{12}{2}
 \mynode{15}{12}{1}
 \mynode{10}{13}{2}
 \mynode{11}{13}{1}
 \mynode{6}{14}{2}
 \mynode{7}{14}{1}
 \mynode{2}{15}{2}
 \mynode{3}{15}{1}
 \mynode{-2}{16}{1}
 \mynode{-1}{16}{1}
\mygrid{-2}{ 12}{ 15}{ 16}
\mylabel{-2}{ 12}{ 15}{ 16}{ 9_{7}^{1}}
\end{tikzpicture}
$
\\
$
\begin{tikzpicture}[yscale=.25,xscale=.25]
\mynode{0}{3}{1}
 \mynode{1}{3}{2}
 \mynode{2}{3}{1}
 \mynode{-4}{4}{1}
 \mynode{-3}{4}{3}
 \mynode{-2}{4}{3}
 \mynode{-1}{4}{1}
 \mynode{-7}{5}{1}
 \mynode{-6}{5}{1}
\mygrid{-7}{ 3}{ 2}{ 5}
\mylabel{-7}{ 3}{ 2}{ 5}{ 9_{7}^{2}}
\end{tikzpicture}
$
\\
$
\begin{tikzpicture}[yscale=.25,xscale=.25]
\mynode{10}{8}{2}
 \mynode{11}{8}{2}
 \mynode{12}{8}{1}
 \mynode{6}{9}{2}
 \mynode{7}{9}{3}
 \mynode{8}{9}{2}
 \mynode{2}{10}{1}
 \mynode{3}{10}{1}
 \mynode{4}{10}{1}
\mygrid{2}{ 8}{ 12}{ 10}
\mylabel{2}{ 8}{ 12}{ 10}{ 9_{8}^{1}}
\end{tikzpicture}
$
\\
$
\begin{tikzpicture}[yscale=.25,xscale=.25]
\mynode{-3}{-13}{1}
 \mynode{-7}{-12}{2}
 \mynode{-6}{-12}{1}
 \mynode{-11}{-11}{2}
 \mynode{-10}{-11}{1}
 \mynode{-15}{-10}{2}
 \mynode{-14}{-10}{1}
 \mynode{-19}{-9}{2}
 \mynode{-18}{-9}{1}
 \mynode{-23}{-8}{1}
 \mynode{-22}{-8}{1}
\mygrid{-23}{ -13}{ -3}{ -8}
\mylabel{-23}{ -13}{ -3}{ -8}{ 9_{8}^{2}}
\end{tikzpicture}
$
\\
$
\begin{tikzpicture}[yscale=.25,xscale=.25]
\mynode{5}{7}{1}
 \mynode{6}{7}{2}
 \mynode{7}{7}{1}
 \mynode{1}{8}{1}
 \mynode{2}{8}{4}
 \mynode{3}{8}{3}
 \mynode{4}{8}{1}
 \mynode{-2}{9}{1}
 \mynode{-1}{9}{1}
\mygrid{-2}{ 7}{ 7}{ 9}
\mylabel{-2}{ 7}{ 7}{ 9}{ 9_{9}^{1}}
\end{tikzpicture}
$
\\
$
\begin{tikzpicture}[yscale=.25,xscale=.25]
\mynode{2}{8}{2}
 \mynode{3}{8}{2}
 \mynode{4}{8}{1}
 \mynode{-2}{9}{2}
 \mynode{-1}{9}{3}
 \mynode{0}{9}{1}
 \mynode{-6}{10}{1}
 \mynode{-5}{10}{2}
 \mynode{-4}{10}{1}
\mygrid{-6}{ 8}{ 4}{ 10}
\mylabel{-6}{ 8}{ 4}{ 10}{ 9_{9}^{2}}
\end{tikzpicture}
$
\\
$
\begin{tikzpicture}[yscale=.25,xscale=.25]
\mynode{0}{-11}{1}
 \mynode{1}{-11}{2}
 \mynode{2}{-11}{1}
 \mynode{-4}{-10}{1}
 \mynode{-3}{-10}{3}
 \mynode{-2}{-10}{1}
 \mynode{-8}{-9}{1}
 \mynode{-7}{-9}{3}
 \mynode{-6}{-9}{1}
 \mynode{-11}{-8}{1}
 \mynode{-10}{-8}{1}
\mygrid{-11}{ -11}{ 2}{ -8}
\mylabel{-11}{ -11}{ 2}{ -8}{ 9_{10}^{1}}
\end{tikzpicture}
$
\\
$
\begin{tikzpicture}[yscale=.25,xscale=.25]
\mynode{0}{-11}{1}
 \mynode{1}{-11}{1}
 \mynode{-4}{-10}{1}
 \mynode{-3}{-10}{2}
 \mynode{-8}{-9}{1}
 \mynode{-7}{-9}{3}
 \mynode{-6}{-9}{2}
 \mynode{-5}{-9}{1}
 \mynode{-11}{-8}{1}
 \mynode{-10}{-8}{2}
 \mynode{-9}{-8}{1}
\mygrid{-11}{ -11}{ 1}{ -8}
\mylabel{-11}{ -11}{ 1}{ -8}{ 9_{10}^{2}}
\end{tikzpicture}
$
\\
$
\begin{tikzpicture}[yscale=.25,xscale=.25]
\mynode{8}{-1}{1}
 \mynode{9}{-1}{2}
 \mynode{10}{-1}{1}
 \mynode{4}{0}{1}
 \mynode{5}{0}{2}
 \mynode{6}{0}{2}
 \mynode{0}{1}{1}
 \mynode{1}{1}{2}
 \mynode{2}{1}{2}
 \mynode{-3}{2}{1}
 \mynode{-2}{2}{1}
\mygrid{-3}{ -1}{ 10}{ 2}
\mylabel{-3}{ -1}{ 10}{ 2}{ 9_{11}^{1}}
\end{tikzpicture}
$
\\
$
\begin{tikzpicture}[yscale=.25,xscale=.25]
\mynode{-1}{-11}{1}
 \mynode{0}{-11}{2}
 \mynode{1}{-11}{1}
 \mynode{-5}{-10}{2}
 \mynode{-4}{-10}{3}
 \mynode{-3}{-10}{2}
 \mynode{-9}{-9}{1}
 \mynode{-8}{-9}{2}
 \mynode{-7}{-9}{2}
\mygrid{-9}{ -11}{ 1}{ -9}
\mylabel{-9}{ -11}{ 1}{ -9}{ 9_{11}^{2}}
\end{tikzpicture}
$
\\
$
\begin{tikzpicture}[yscale=.25,xscale=.25]
\mynode{2}{4}{2}
 \mynode{3}{4}{1}
 \mynode{-2}{5}{2}
 \mynode{-1}{5}{2}
 \mynode{-6}{6}{2}
 \mynode{-5}{6}{2}
 \mynode{-10}{7}{2}
 \mynode{-9}{7}{2}
 \mynode{-14}{8}{1}
 \mynode{-13}{8}{1}
\mygrid{-14}{ 4}{ 3}{ 8}
\mylabel{-14}{ 4}{ 3}{ 8}{ 9_{12}^{1}}
\end{tikzpicture}
$
\\
$
\begin{tikzpicture}[yscale=.25,xscale=.25]
\mynode{12}{1}{1}
 \mynode{13}{1}{1}
 \mynode{8}{2}{2}
 \mynode{9}{2}{3}
 \mynode{10}{2}{1}
 \mynode{4}{3}{2}
 \mynode{5}{3}{3}
 \mynode{6}{3}{1}
 \mynode{0}{4}{1}
 \mynode{1}{4}{1}
 \mynode{2}{4}{1}
\mygrid{0}{ 1}{ 13}{ 4}
\mylabel{0}{ 1}{ 13}{ 4}{ 9_{12}^{2}}
\end{tikzpicture}
$
\\
$
\begin{tikzpicture}[yscale=.25,xscale=.25]
\mynode{8}{-7}{1}
 \mynode{9}{-7}{2}
 \mynode{10}{-7}{1}
 \mynode{4}{-6}{1}
 \mynode{5}{-6}{3}
 \mynode{6}{-6}{2}
 \mynode{0}{-5}{1}
 \mynode{1}{-5}{3}
 \mynode{2}{-5}{2}
 \mynode{-3}{-4}{1}
 \mynode{-2}{-4}{1}
\mygrid{-3}{ -7}{ 10}{ -4}
\mylabel{-3}{ -7}{ 10}{ -4}{ 9_{13}^{1}}
\end{tikzpicture}
$
\\
$
\begin{tikzpicture}[yscale=.25,xscale=.25]
\mynode{-3}{-16}{1}
 \mynode{-8}{-15}{1}
 \mynode{-7}{-15}{3}
 \mynode{-6}{-15}{1}
 \mynode{-12}{-14}{1}
 \mynode{-11}{-14}{3}
 \mynode{-10}{-14}{1}
 \mynode{-16}{-13}{1}
 \mynode{-15}{-13}{3}
 \mynode{-14}{-13}{1}
 \mynode{-19}{-12}{1}
 \mynode{-18}{-12}{1}
\mygrid{-19}{ -16}{ -3}{ -12}
\mylabel{-19}{ -16}{ -3}{ -12}{ 9_{13}^{2}}
\end{tikzpicture}
$
\\
$
\begin{tikzpicture}[yscale=.25,xscale=.25]
\mynode{21}{0}{1}
 \mynode{22}{0}{1}
 \mynode{17}{1}{2}
 \mynode{18}{1}{2}
 \mynode{13}{2}{2}
 \mynode{14}{2}{2}
 \mynode{9}{3}{2}
 \mynode{10}{3}{2}
 \mynode{5}{4}{1}
 \mynode{6}{4}{2}
 \mynode{2}{5}{1}
\mygrid{2}{ 0}{ 22}{ 5}
\mylabel{2}{ 0}{ 22}{ 5}{ 9_{14}^{1}}
\end{tikzpicture}
$
\\
$
\begin{tikzpicture}[yscale=.25,xscale=.25]
\mynode{-3}{-6}{1}
 \mynode{-2}{-6}{1}
 \mynode{-1}{-6}{1}
 \mynode{-7}{-5}{2}
 \mynode{-6}{-5}{3}
 \mynode{-5}{-5}{2}
 \mynode{-11}{-4}{1}
 \mynode{-10}{-4}{3}
 \mynode{-9}{-4}{2}
 \mynode{-14}{-3}{1}
 \mynode{-13}{-3}{1}
\mygrid{-14}{ -6}{ -1}{ -3}
\mylabel{-14}{ -6}{ -1}{ -3}{ 9_{14}^{2}}
\end{tikzpicture}
$
\\
$
\begin{tikzpicture}[yscale=.25,xscale=.25]
\mynode{5}{-1}{1}
 \mynode{6}{-1}{2}
 \mynode{7}{-1}{1}
 \mynode{1}{0}{1}
 \mynode{2}{0}{4}
 \mynode{3}{0}{2}
 \mynode{-3}{1}{1}
 \mynode{-2}{1}{3}
 \mynode{-1}{1}{2}
 \mynode{-6}{2}{1}
 \mynode{-5}{2}{1}
\mygrid{-6}{ -1}{ 7}{ 2}
\mylabel{-6}{ -1}{ 7}{ 2}{ 9_{15}^{1}}
\end{tikzpicture}
$
\\
$
\begin{tikzpicture}[yscale=.25,xscale=.25]
\mynode{4}{-13}{1}
 \mynode{5}{-13}{1}
 \mynode{0}{-12}{2}
 \mynode{1}{-12}{2}
 \mynode{-4}{-11}{2}
 \mynode{-3}{-11}{3}
 \mynode{-8}{-10}{2}
 \mynode{-7}{-10}{3}
 \mynode{-12}{-9}{1}
 \mynode{-11}{-9}{2}
\mygrid{-12}{ -13}{ 5}{ -9}
\mylabel{-12}{ -13}{ 5}{ -9}{ 9_{15}^{2}}
\end{tikzpicture}
$
\\
$
\begin{tikzpicture}[yscale=.25,xscale=.25]
\mynode{3}{-11}{1}
 \mynode{4}{-11}{2}
 \mynode{5}{-11}{1}
 \mynode{-1}{-10}{1}
 \mynode{0}{-10}{4}
 \mynode{1}{-10}{3}
 \mynode{-5}{-9}{1}
 \mynode{-4}{-9}{3}
 \mynode{-3}{-9}{3}
\mygrid{-5}{ -11}{ 5}{ -9}
\mylabel{-5}{ -11}{ 5}{ -9}{ 9_{16}^{1}}
\end{tikzpicture}
$
\\
$
\begin{tikzpicture}[yscale=.25,xscale=.25]
\mynode{16}{-2}{1}
 \mynode{17}{-2}{1}
 \mynode{11}{-1}{1}
 \mynode{12}{-1}{3}
 \mynode{13}{-1}{3}
 \mynode{14}{-1}{1}
 \mynode{8}{0}{1}
 \mynode{9}{0}{2}
 \mynode{10}{0}{2}
 \mynode{5}{1}{1}
 \mynode{6}{1}{2}
 \mynode{2}{2}{1}
\mygrid{2}{ -2}{ 17}{ 2}
\mylabel{2}{ -2}{ 17}{ 2}{ 9_{16}^{2}}
\end{tikzpicture}
$
\\
$
\begin{tikzpicture}[yscale=.25,xscale=.25]
\mynode{1}{-2}{1}
 \mynode{2}{-2}{2}
 \mynode{3}{-2}{1}
 \mynode{-3}{-1}{2}
 \mynode{-2}{-1}{3}
 \mynode{-1}{-1}{2}
 \mynode{-7}{0}{1}
 \mynode{-6}{0}{3}
 \mynode{-5}{0}{2}
 \mynode{-10}{1}{1}
 \mynode{-9}{1}{1}
\mygrid{-10}{ -2}{ 3}{ 1}
\mylabel{-10}{ -2}{ 3}{ 1}{ 9_{17}^{1}}
\end{tikzpicture}
$
\\
$
\begin{tikzpicture}[yscale=.25,xscale=.25]
\mynode{1}{-2}{1}
 \mynode{2}{-2}{1}
 \mynode{-3}{-1}{2}
 \mynode{-2}{-1}{2}
 \mynode{-1}{-1}{1}
 \mynode{-7}{0}{1}
 \mynode{-6}{0}{3}
 \mynode{-5}{0}{3}
 \mynode{-4}{0}{1}
 \mynode{-10}{1}{1}
 \mynode{-9}{1}{2}
 \mynode{-8}{1}{1}
\mygrid{-10}{ -2}{ 2}{ 1}
\mylabel{-10}{ -2}{ 2}{ 1}{ 9_{17}^{2}}
\end{tikzpicture}
$
\\
$
\begin{tikzpicture}[yscale=.25,xscale=.25]
\mynode{6}{8}{2}
 \mynode{7}{8}{1}
 \mynode{2}{9}{3}
 \mynode{3}{9}{2}
 \mynode{-2}{10}{3}
 \mynode{-1}{10}{3}
 \mynode{-6}{11}{2}
 \mynode{-5}{11}{2}
 \mynode{-10}{12}{1}
 \mynode{-9}{12}{1}
\mygrid{-10}{ 8}{ 7}{ 12}
\mylabel{-10}{ 8}{ 7}{ 12}{ 9_{18}^{1}}
\end{tikzpicture}
$
\\
$
\begin{tikzpicture}[yscale=.25,xscale=.25]
\mynode{-3}{-2}{1}
 \mynode{-8}{-1}{1}
 \mynode{-7}{-1}{4}
 \mynode{-6}{-1}{2}
 \mynode{-12}{0}{1}
 \mynode{-11}{0}{4}
 \mynode{-10}{0}{4}
 \mynode{-9}{0}{1}
 \mynode{-15}{1}{1}
 \mynode{-14}{1}{1}
\mygrid{-15}{ -2}{ -3}{ 1}
\mylabel{-15}{ -2}{ -3}{ 1}{ 9_{18}^{2}}
\end{tikzpicture}
$
\\
$
\begin{tikzpicture}[yscale=.25,xscale=.25]
\mynode{4}{-2}{1}
 \mynode{5}{-2}{2}
 \mynode{6}{-2}{1}
 \mynode{0}{-1}{2}
 \mynode{1}{-1}{4}
 \mynode{2}{-1}{2}
 \mynode{-4}{0}{1}
 \mynode{-3}{0}{3}
 \mynode{-2}{0}{2}
 \mynode{-7}{1}{1}
 \mynode{-6}{1}{1}
\mygrid{-7}{ -2}{ 6}{ 1}
\mylabel{-7}{ -2}{ 6}{ 1}{ 9_{19}^{1}}
\end{tikzpicture}
$
\\
$
\begin{tikzpicture}[yscale=.25,xscale=.25]
\mynode{21}{6}{1}
 \mynode{22}{6}{1}
 \mynode{17}{7}{2}
 \mynode{18}{7}{2}
 \mynode{13}{8}{2}
 \mynode{14}{8}{3}
 \mynode{9}{9}{2}
 \mynode{10}{9}{3}
 \mynode{5}{10}{1}
 \mynode{6}{10}{2}
 \mynode{2}{11}{1}
\mygrid{2}{ 6}{ 22}{ 11}
\mylabel{2}{ 6}{ 22}{ 11}{ 9_{19}^{2}}
\end{tikzpicture}
$
\\
$
\begin{tikzpicture}[yscale=.25,xscale=.25]
\mynode{6}{8}{3}
 \mynode{7}{8}{3}
 \mynode{8}{8}{1}
 \mynode{2}{9}{3}
 \mynode{3}{9}{4}
 \mynode{4}{9}{2}
 \mynode{-2}{10}{1}
 \mynode{-1}{10}{2}
 \mynode{0}{10}{1}
\mygrid{-2}{ 8}{ 8}{ 10}
\mylabel{-2}{ 8}{ 8}{ 10}{ 9_{20}^{1}}
\end{tikzpicture}
$
\\
$
\begin{tikzpicture}[yscale=.25,xscale=.25]
\mynode{17}{5}{1}
 \mynode{18}{5}{1}
 \mynode{13}{6}{2}
 \mynode{14}{6}{3}
 \mynode{15}{6}{1}
 \mynode{9}{7}{2}
 \mynode{10}{7}{3}
 \mynode{11}{7}{1}
 \mynode{5}{8}{1}
 \mynode{6}{8}{3}
 \mynode{7}{8}{1}
 \mynode{2}{9}{1}
\mygrid{2}{ 5}{ 18}{ 9}
\mylabel{2}{ 5}{ 18}{ 9}{ 9_{20}^{2}}
\end{tikzpicture}
$
\\
$
\begin{tikzpicture}[yscale=.25,xscale=.25]
\mynode{21}{1}{1}
 \mynode{22}{1}{1}
 \mynode{17}{2}{2}
 \mynode{18}{2}{2}
 \mynode{13}{3}{3}
 \mynode{14}{3}{3}
 \mynode{9}{4}{2}
 \mynode{10}{4}{3}
 \mynode{5}{5}{1}
 \mynode{6}{5}{2}
 \mynode{2}{6}{1}
\mygrid{2}{ 1}{ 22}{ 6}
\mylabel{2}{ 1}{ 22}{ 6}{ 9_{21}^{1}}
\end{tikzpicture}
$
\\
$
\begin{tikzpicture}[yscale=.25,xscale=.25]
\mynode{4}{-7}{1}
 \mynode{5}{-7}{2}
 \mynode{6}{-7}{1}
 \mynode{0}{-6}{2}
 \mynode{1}{-6}{4}
 \mynode{2}{-6}{2}
 \mynode{-4}{-5}{1}
 \mynode{-3}{-5}{4}
 \mynode{-2}{-5}{2}
 \mynode{-7}{-4}{1}
 \mynode{-6}{-4}{1}
\mygrid{-7}{ -7}{ 6}{ -4}
\mylabel{-7}{ -7}{ 6}{ -4}{ 9_{21}^{2}}
\end{tikzpicture}
$
\\
$
\begin{tikzpicture}[yscale=.25,xscale=.25]
\mynode{0}{-6}{1}
 \mynode{1}{-6}{1}
 \mynode{-4}{-5}{2}
 \mynode{-3}{-5}{4}
 \mynode{-2}{-5}{1}
 \mynode{-8}{-4}{2}
 \mynode{-7}{-4}{4}
 \mynode{-6}{-4}{2}
 \mynode{-12}{-3}{1}
 \mynode{-11}{-3}{2}
 \mynode{-10}{-3}{1}
\mygrid{-12}{ -6}{ 1}{ -3}
\mylabel{-12}{ -6}{ 1}{ -3}{ 9_{22}}
\end{tikzpicture}
$
\\
$
\begin{tikzpicture}[yscale=.25,xscale=.25]
\mynode{17}{11}{1}
 \mynode{18}{11}{1}
 \mynode{13}{12}{2}
 \mynode{14}{12}{4}
 \mynode{15}{12}{1}
 \mynode{9}{13}{2}
 \mynode{10}{13}{4}
 \mynode{11}{13}{1}
 \mynode{5}{14}{1}
 \mynode{6}{14}{3}
 \mynode{7}{14}{1}
 \mynode{2}{15}{1}
\mygrid{2}{ 11}{ 18}{ 15}
\mylabel{2}{ 11}{ 18}{ 15}{ 9_{23}^{1}}
\end{tikzpicture}
$
\\
$
\begin{tikzpicture}[yscale=.25,xscale=.25]
\mynode{4}{2}{1}
 \mynode{5}{2}{1}
 \mynode{0}{3}{2}
 \mynode{1}{3}{3}
 \mynode{2}{3}{1}
 \mynode{-4}{4}{1}
 \mynode{-3}{4}{3}
 \mynode{-2}{4}{3}
 \mynode{-1}{4}{1}
 \mynode{-7}{5}{1}
 \mynode{-6}{5}{2}
 \mynode{-5}{5}{1}
 \mynode{-10}{6}{1}
 \mynode{-9}{6}{1}
\mygrid{-10}{ 2}{ 5}{ 6}
\mylabel{-10}{ 2}{ 5}{ 6}{ 9_{23}^{2}}
\end{tikzpicture}
$
\\
$
\begin{tikzpicture}[yscale=.25,xscale=.25]
\mynode{0}{-7}{1}
 \mynode{1}{-7}{1}
 \mynode{-4}{-6}{2}
 \mynode{-3}{-6}{4}
 \mynode{-2}{-6}{1}
 \mynode{-8}{-5}{2}
 \mynode{-7}{-5}{5}
 \mynode{-6}{-5}{2}
 \mynode{-12}{-4}{1}
 \mynode{-11}{-4}{2}
 \mynode{-10}{-4}{1}
\mygrid{-12}{ -7}{ 1}{ -4}
\mylabel{-12}{ -7}{ 1}{ -4}{ 9_{24}}
\end{tikzpicture}
$
\\
$
\begin{tikzpicture}[yscale=.25,xscale=.25]
\mynode{4}{-3}{1}
 \mynode{5}{-3}{1}
 \mynode{0}{-2}{3}
 \mynode{1}{-2}{4}
 \mynode{2}{-2}{1}
 \mynode{-4}{-1}{3}
 \mynode{-3}{-1}{5}
 \mynode{-2}{-1}{1}
 \mynode{-8}{0}{1}
 \mynode{-7}{0}{2}
 \mynode{-6}{0}{1}
\mygrid{-8}{ -3}{ 5}{ 0}
\mylabel{-8}{ -3}{ 5}{ 0}{ 9_{25}}
\end{tikzpicture}
$
\\
$
\begin{tikzpicture}[yscale=.25,xscale=.25]
\mynode{8}{-2}{1}
 \mynode{9}{-2}{2}
 \mynode{10}{-2}{1}
 \mynode{4}{-1}{2}
 \mynode{5}{-1}{5}
 \mynode{6}{-1}{3}
 \mynode{0}{0}{1}
 \mynode{1}{0}{3}
 \mynode{2}{0}{3}
 \mynode{-3}{1}{1}
 \mynode{-2}{1}{1}
\mygrid{-3}{ -2}{ 10}{ 1}
\mylabel{-3}{ -2}{ 10}{ 1}{ 9_{26}^{1}}
\end{tikzpicture}
$
\\
$
\begin{tikzpicture}[yscale=.25,xscale=.25]
\mynode{-3}{-11}{1}
 \mynode{-8}{-10}{1}
 \mynode{-7}{-10}{3}
 \mynode{-6}{-10}{1}
 \mynode{-12}{-9}{2}
 \mynode{-11}{-9}{5}
 \mynode{-10}{-9}{2}
 \mynode{-16}{-8}{1}
 \mynode{-15}{-8}{3}
 \mynode{-14}{-8}{2}
 \mynode{-19}{-7}{1}
 \mynode{-18}{-7}{1}
\mygrid{-19}{ -11}{ -3}{ -7}
\mylabel{-19}{ -11}{ -3}{ -7}{ 9_{26}^{2}}
\end{tikzpicture}
$
\\
$
\begin{tikzpicture}[yscale=.25,xscale=.25]
\mynode{17}{7}{1}
 \mynode{18}{7}{1}
 \mynode{13}{8}{2}
 \mynode{14}{8}{4}
 \mynode{15}{8}{1}
 \mynode{9}{9}{2}
 \mynode{10}{9}{5}
 \mynode{11}{9}{2}
 \mynode{5}{10}{1}
 \mynode{6}{10}{3}
 \mynode{7}{10}{1}
 \mynode{2}{11}{1}
\mygrid{2}{ 7}{ 18}{ 11}
\mylabel{2}{ 7}{ 18}{ 11}{ 9_{27}^{1}}
\end{tikzpicture}
$
\\
$
\begin{tikzpicture}[yscale=.25,xscale=.25]
\mynode{8}{-3}{1}
 \mynode{9}{-3}{1}
 \mynode{4}{-2}{3}
 \mynode{5}{-2}{4}
 \mynode{6}{-2}{1}
 \mynode{0}{-1}{3}
 \mynode{1}{-1}{5}
 \mynode{2}{-1}{2}
 \mynode{-4}{0}{1}
 \mynode{-3}{0}{2}
 \mynode{-2}{0}{1}
\mygrid{-4}{ -3}{ 9}{ 0}
\mylabel{-4}{ -3}{ 9}{ 0}{ 9_{27}^{2}}
\end{tikzpicture}
$
\\
$
\begin{tikzpicture}[yscale=.25,xscale=.25]
\mynode{-3}{-7}{1}
 \mynode{-7}{-6}{3}
 \mynode{-6}{-6}{2}
 \mynode{-11}{-5}{4}
 \mynode{-10}{-5}{4}
 \mynode{-9}{-5}{1}
 \mynode{-15}{-4}{3}
 \mynode{-14}{-4}{4}
 \mynode{-13}{-4}{1}
 \mynode{-19}{-3}{1}
 \mynode{-18}{-3}{1}
\mygrid{-19}{ -7}{ -3}{ -3}
\mylabel{-19}{ -7}{ -3}{ -3}{ 9_{28}^{1}}
\end{tikzpicture}
$
\\
$
\begin{tikzpicture}[yscale=.25,xscale=.25]
\mynode{-3}{-7}{1}
 \mynode{-8}{-6}{1}
 \mynode{-7}{-6}{3}
 \mynode{-6}{-6}{1}
 \mynode{-12}{-5}{1}
 \mynode{-11}{-5}{4}
 \mynode{-10}{-5}{3}
 \mynode{-9}{-5}{1}
 \mynode{-15}{-4}{2}
 \mynode{-14}{-4}{4}
 \mynode{-13}{-4}{2}
 \mynode{-18}{-3}{1}
 \mynode{-17}{-3}{1}
\mygrid{-18}{ -7}{ -3}{ -3}
\mylabel{-18}{ -7}{ -3}{ -3}{ 9_{28}^{2}}
\end{tikzpicture}
$
\\
$
\begin{tikzpicture}[yscale=.25,xscale=.25]
\mynode{-3}{-7}{1}
 \mynode{-7}{-6}{2}
 \mynode{-6}{-6}{1}
 \mynode{-11}{-5}{2}
 \mynode{-10}{-5}{3}
 \mynode{-9}{-5}{2}
 \mynode{-15}{-4}{1}
 \mynode{-14}{-4}{4}
 \mynode{-13}{-4}{4}
 \mynode{-12}{-4}{1}
 \mynode{-18}{-3}{1}
 \mynode{-17}{-3}{2}
 \mynode{-16}{-3}{1}
\mygrid{-18}{ -7}{ -3}{ -3}
\mylabel{-18}{ -7}{ -3}{ -3}{ 9_{29}}
\end{tikzpicture}
$
\\
$
\begin{tikzpicture}[yscale=.25,xscale=.25]
\mynode{-3}{-6}{1}
 \mynode{-7}{-5}{3}
 \mynode{-6}{-5}{2}
 \mynode{-11}{-4}{4}
 \mynode{-10}{-4}{5}
 \mynode{-9}{-4}{1}
 \mynode{-15}{-3}{3}
 \mynode{-14}{-3}{4}
 \mynode{-13}{-3}{1}
 \mynode{-19}{-2}{1}
 \mynode{-18}{-2}{1}
\mygrid{-19}{ -6}{ -3}{ -2}
\mylabel{-19}{ -6}{ -3}{ -2}{ 9_{30}}
\end{tikzpicture}
$
\\
$
\begin{tikzpicture}[yscale=.25,xscale=.25]
\mynode{-3}{-11}{1}
 \mynode{-8}{-10}{1}
 \mynode{-7}{-10}{3}
 \mynode{-6}{-10}{1}
 \mynode{-12}{-9}{2}
 \mynode{-11}{-9}{6}
 \mynode{-10}{-9}{3}
 \mynode{-16}{-8}{1}
 \mynode{-15}{-8}{4}
 \mynode{-14}{-8}{3}
 \mynode{-19}{-7}{1}
 \mynode{-18}{-7}{1}
\mygrid{-19}{ -11}{ -3}{ -7}
\mylabel{-19}{ -11}{ -3}{ -7}{ 9_{31}^{1}}
\end{tikzpicture}
$
\\
$
\begin{tikzpicture}[yscale=.25,xscale=.25]
\mynode{8}{-2}{1}
 \mynode{9}{-2}{1}
 \mynode{4}{-1}{2}
 \mynode{5}{-1}{4}
 \mynode{6}{-1}{2}
 \mynode{0}{0}{1}
 \mynode{1}{0}{3}
 \mynode{2}{0}{4}
 \mynode{3}{0}{1}
 \mynode{-3}{1}{1}
 \mynode{-2}{1}{3}
 \mynode{-1}{1}{2}
 \mynode{-6}{2}{1}
 \mynode{-5}{2}{1}
\mygrid{-6}{ -2}{ 9}{ 2}
\mylabel{-6}{ -2}{ 9}{ 2}{ 9_{31}^{2}}
\end{tikzpicture}
$
\\
$
\begin{tikzpicture}[yscale=.25,xscale=.25]
\mynode{-3}{-6}{1}
 \mynode{-7}{-5}{3}
 \mynode{-6}{-5}{2}
 \mynode{-11}{-4}{4}
 \mynode{-10}{-4}{5}
 \mynode{-9}{-4}{1}
 \mynode{-15}{-3}{3}
 \mynode{-14}{-3}{5}
 \mynode{-13}{-3}{2}
 \mynode{-19}{-2}{1}
 \mynode{-18}{-2}{3}
 \mynode{-17}{-2}{2}
 \mynode{-22}{-1}{1}
 \mynode{-21}{-1}{1}
\mygrid{-22}{ -6}{ -3}{ -1}
\mylabel{-22}{ -6}{ -3}{ -1}{ 9_{34}}
\end{tikzpicture}
$
\\
$
\begin{tikzpicture}[yscale=.25,xscale=.25]
\mynode{20}{11}{1}
 \mynode{21}{11}{1}
 \mynode{16}{12}{1}
 \mynode{17}{12}{1}
 \mynode{18}{12}{1}
 \mynode{12}{13}{1}
 \mynode{13}{13}{1}
 \mynode{14}{13}{1}
 \mynode{8}{14}{1}
 \mynode{9}{14}{1}
 \mynode{10}{14}{1}
 \mynode{6}{15}{1}
 \mynode{2}{16}{1}
\mygrid{2}{ 11}{ 21}{ 16}
\mylabel{2}{ 11}{ 21}{ 16}{ 9_{35}^{1}}
\end{tikzpicture}
$
\\
$
\begin{tikzpicture}[yscale=.25,xscale=.25]
\mynode{3}{3}{1}
 \mynode{4}{3}{2}
 \mynode{5}{3}{1}
 \mynode{-1}{4}{1}
 \mynode{0}{4}{2}
 \mynode{1}{4}{1}
 \mynode{2}{4}{1}
 \mynode{-3}{5}{1}
 \mynode{-2}{5}{1}
 \mynode{-7}{6}{1}
 \mynode{-6}{6}{1}
\mygrid{-7}{ 3}{ 5}{ 6}
\mylabel{-7}{ 3}{ 5}{ 6}{ 9_{35}^{2}}
\end{tikzpicture}
$
\\
$
\begin{tikzpicture}[yscale=.25,xscale=.25]
\mynode{7}{-7}{1}
 \mynode{8}{-7}{2}
 \mynode{9}{-7}{1}
 \mynode{3}{-6}{2}
 \mynode{4}{-6}{4}
 \mynode{5}{-6}{2}
 \mynode{-1}{-5}{1}
 \mynode{0}{-5}{3}
 \mynode{1}{-5}{2}
\mygrid{-1}{ -7}{ 9}{ -5}
\mylabel{-1}{ -7}{ 9}{ -5}{ 9_{36}}
\end{tikzpicture}
$
\\
$
\begin{tikzpicture}[yscale=.25,xscale=.25]
\mynode{-3}{-6}{1}
 \mynode{-7}{-5}{2}
 \mynode{-6}{-5}{1}
 \mynode{-11}{-4}{2}
 \mynode{-10}{-4}{3}
 \mynode{-9}{-4}{1}
 \mynode{-15}{-3}{2}
 \mynode{-14}{-3}{3}
 \mynode{-13}{-3}{1}
 \mynode{-19}{-2}{1}
 \mynode{-18}{-2}{2}
 \mynode{-17}{-2}{1}
 \mynode{-22}{-1}{1}
 \mynode{-21}{-1}{1}
\mygrid{-22}{ -6}{ -3}{ -1}
\mylabel{-22}{ -6}{ -3}{ -1}{ 9_{37}^{1}}
\end{tikzpicture}
$
\\
$
\begin{tikzpicture}[yscale=.25,xscale=.25]
\mynode{-4}{-6}{1}
 \mynode{-3}{-6}{2}
 \mynode{-8}{-5}{2}
 \mynode{-7}{-5}{4}
 \mynode{-6}{-5}{1}
 \mynode{-12}{-4}{1}
 \mynode{-11}{-4}{3}
 \mynode{-10}{-4}{3}
 \mynode{-9}{-4}{1}
 \mynode{-15}{-3}{1}
 \mynode{-14}{-3}{2}
 \mynode{-13}{-3}{1}
\mygrid{-15}{ -6}{ -3}{ -3}
\mylabel{-15}{ -6}{ -3}{ -3}{ 9_{37}^{2}}
\end{tikzpicture}
$
\\
$
\begin{tikzpicture}[yscale=.25,xscale=.25]
\mynode{12}{6}{1}
 \mynode{13}{6}{1}
 \mynode{8}{7}{2}
 \mynode{9}{7}{4}
 \mynode{10}{7}{2}
 \mynode{4}{8}{1}
 \mynode{5}{8}{4}
 \mynode{6}{8}{5}
 \mynode{7}{8}{1}
 \mynode{1}{9}{1}
 \mynode{2}{9}{3}
 \mynode{3}{9}{1}
 \mynode{-2}{10}{1}
 \mynode{-1}{10}{1}
\mygrid{-2}{ 6}{ 13}{ 10}
\mylabel{-2}{ 6}{ 13}{ 10}{ 9_{38}}
\end{tikzpicture}
$
\\
$
\begin{tikzpicture}[yscale=.25,xscale=.25]
\mynode{-3}{-11}{1}
 \mynode{-7}{-10}{2}
 \mynode{-6}{-10}{1}
 \mynode{-11}{-9}{3}
 \mynode{-10}{-9}{3}
 \mynode{-9}{-9}{1}
 \mynode{-15}{-8}{2}
 \mynode{-14}{-8}{4}
 \mynode{-13}{-8}{2}
 \mynode{-19}{-7}{1}
 \mynode{-18}{-7}{3}
 \mynode{-17}{-7}{2}
 \mynode{-22}{-6}{1}
 \mynode{-21}{-6}{1}
\mygrid{-22}{ -11}{ -3}{ -6}
\mylabel{-22}{ -11}{ -3}{ -6}{ 9_{39}}
\end{tikzpicture}
$
\\
$
\begin{tikzpicture}[yscale=.25,xscale=.25]
\mynode{-3}{-7}{1}
 \mynode{-7}{-6}{3}
 \mynode{-6}{-6}{2}
 \mynode{-11}{-5}{4}
 \mynode{-10}{-5}{4}
 \mynode{-9}{-5}{1}
 \mynode{-15}{-4}{3}
 \mynode{-14}{-4}{6}
 \mynode{-13}{-4}{3}
 \mynode{-19}{-3}{1}
 \mynode{-18}{-3}{4}
 \mynode{-17}{-3}{3}
 \mynode{-22}{-2}{1}
 \mynode{-21}{-2}{1}
\mygrid{-22}{ -7}{ -3}{ -2}
\mylabel{-22}{ -7}{ -3}{ -2}{ 9_{40}^{1}}
\end{tikzpicture}
$
\\
$
\begin{tikzpicture}[yscale=.25,xscale=.25]
\mynode{12}{1}{1}
 \mynode{13}{1}{1}
 \mynode{8}{2}{2}
 \mynode{9}{2}{4}
 \mynode{10}{2}{2}
 \mynode{4}{3}{1}
 \mynode{5}{3}{5}
 \mynode{6}{3}{5}
 \mynode{7}{3}{1}
 \mynode{1}{4}{2}
 \mynode{2}{4}{6}
 \mynode{3}{4}{3}
 \mynode{-2}{5}{2}
 \mynode{-1}{5}{2}
\mygrid{-2}{ 1}{ 13}{ 5}
\mylabel{-2}{ 1}{ 13}{ 5}{ 9_{40}^{2}}
\end{tikzpicture}
$
\\
$
\begin{tikzpicture}[yscale=.25,xscale=.25]
\mynode{20}{6}{1}
 \mynode{21}{6}{1}
 \mynode{16}{7}{2}
 \mynode{17}{7}{3}
 \mynode{18}{7}{1}
 \mynode{12}{8}{2}
 \mynode{13}{8}{3}
 \mynode{14}{8}{2}
 \mynode{8}{9}{1}
 \mynode{9}{9}{2}
 \mynode{10}{9}{2}
 \mynode{5}{10}{1}
 \mynode{6}{10}{2}
 \mynode{2}{11}{1}
\mygrid{2}{ 6}{ 21}{ 11}
\mylabel{2}{ 6}{ 21}{ 11}{ 9_{41}}
\end{tikzpicture}
$
\\
$
\begin{tikzpicture}[yscale=.25,xscale=.25]
\mynode{7}{-2}{1}
 \mynode{8}{-2}{1}
 \mynode{3}{-1}{1}
 \mynode{4}{-1}{1}
\mygrid{3}{ -2}{ 8}{ -1}
\mylabel{3}{ -2}{ 8}{ -1}{ 9_{42}}
\end{tikzpicture}
$
\\
$
\begin{tikzpicture}[yscale=.25,xscale=.25]
\mynode{0}{-5}{1}
 \mynode{1}{-5}{1}
 \mynode{-4}{-4}{1}
 \mynode{-3}{-4}{1}
 \mynode{-8}{-3}{1}
 \mynode{-7}{-3}{1}
\mygrid{-8}{ -5}{ 1}{ -3}
\mylabel{-8}{ -5}{ 1}{ -3}{ 9_{43}}
\end{tikzpicture}
$
\\
$
\begin{tikzpicture}[yscale=.25,xscale=.25]
\mynode{4}{-2}{1}
 \mynode{5}{-2}{1}
 \mynode{0}{-1}{2}
 \mynode{1}{-1}{2}
 \mynode{-4}{0}{1}
 \mynode{-3}{0}{1}
\mygrid{-4}{ -2}{ 5}{ 0}
\mylabel{-4}{ -2}{ 5}{ 0}{ 9_{44}}
\end{tikzpicture}
$
\\
$
\begin{tikzpicture}[yscale=.25,xscale=.25]
\mynode{-3}{-1}{1}
 \mynode{-7}{0}{2}
 \mynode{-6}{0}{2}
 \mynode{-11}{1}{2}
 \mynode{-10}{1}{2}
 \mynode{-15}{2}{1}
 \mynode{-14}{2}{1}
\mygrid{-15}{ -1}{ -3}{ 2}
\mylabel{-15}{ -1}{ -3}{ 2}{ 9_{45}}
\end{tikzpicture}
$
\\
$
\begin{tikzpicture}[yscale=.25,xscale=.25]
\mynode{3}{-2}{1}
 \mynode{4}{-2}{1}
 \mynode{-1}{-1}{1}
 \mynode{0}{-1}{1}
\mygrid{-1}{ -2}{ 4}{ -1}
\mylabel{-1}{ -2}{ 4}{ -1}{ 9_{46}^{1}}
\end{tikzpicture}
$
\\
$
\begin{tikzpicture}[yscale=.25,xscale=.25]
\mynode{10}{8}{1}
 \mynode{12}{8}{1}
 \mynode{13}{8}{1}
 \mynode{8}{9}{1}
 \mynode{9}{9}{1}
 \mynode{6}{10}{1}
 \mynode{2}{11}{1}
\mygrid{2}{ 8}{ 13}{ 11}
\mylabel{2}{ 8}{ 13}{ 11}{ 9_{46}^{2}}
\end{tikzpicture}
$
\\
$
\begin{tikzpicture}[yscale=.25,xscale=.25]
\mynode{0}{-6}{1}
 \mynode{1}{-6}{1}
 \mynode{-4}{-5}{2}
 \mynode{-3}{-5}{3}
 \mynode{-8}{-4}{1}
 \mynode{-7}{-4}{2}
 \mynode{-6}{-4}{1}
 \mynode{-11}{-3}{1}
 \mynode{-10}{-3}{1}
\mygrid{-11}{ -6}{ 1}{ -3}
\mylabel{-11}{ -6}{ 1}{ -3}{ 9_{47}}
\end{tikzpicture}
$
\\
$
\begin{tikzpicture}[yscale=.25,xscale=.25]
\mynode{5}{3}{1}
 \mynode{6}{3}{1}
 \mynode{1}{4}{1}
 \mynode{2}{4}{3}
 \mynode{3}{4}{1}
 \mynode{-2}{5}{2}
 \mynode{-1}{5}{2}
 \mynode{-6}{6}{1}
 \mynode{-5}{6}{1}
\mygrid{-6}{ 3}{ 6}{ 6}
\mylabel{-6}{ 3}{ 6}{ 6}{ 9_{48}^{1}}
\end{tikzpicture}
$
\\
$
\begin{tikzpicture}[yscale=.25,xscale=.25]
\mynode{12}{1}{1}
 \mynode{13}{1}{1}
 \mynode{8}{2}{1}
 \mynode{9}{2}{2}
 \mynode{10}{2}{1}
 \mynode{5}{3}{2}
 \mynode{6}{3}{2}
 \mynode{1}{4}{1}
 \mynode{2}{4}{2}
\mygrid{1}{ 1}{ 13}{ 4}
\mylabel{1}{ 1}{ 13}{ 4}{ 9_{48}^{2}}
\end{tikzpicture}
$
\\
$
\begin{tikzpicture}[yscale=.25,xscale=.25]
\mynode{0}{-11}{1}
 \mynode{1}{-11}{1}
 \mynode{-4}{-10}{1}
 \mynode{-3}{-10}{2}
 \mynode{-8}{-9}{1}
 \mynode{-7}{-9}{3}
 \mynode{-6}{-9}{1}
 \mynode{-11}{-8}{1}
 \mynode{-10}{-8}{1}
\mygrid{-11}{ -11}{ 1}{ -8}
\mylabel{-11}{ -11}{ 1}{ -8}{ 9_{49}}
\end{tikzpicture}
$
\\